\input amstex
\magnification=\magstep1
\hsize=16truecm
\parskip=4pt

\define\({\left(}
\define\){\right)}
\define\[{\left[}
\define\]{\right]}
\define\e{\varepsilon}
\define\oo{\omega}

\define\supp {\sup\limits}

\define\summ{\sum\limits}

\define\bigcapp{\bigcap\limits}

\TagsOnRight
 
\centerline{\bf Estimates on the tail behaviour of Gaussian polynomials.}
\centerline {\bf The discussion of a result of Lata{\l}a.}
 
\medskip
\centerline{\it P\'eter Major}
\centerline{Alfr\'ed R\'enyi Mathematical Institute of the Hungarian
Academy of Sciences}
\centerline{Budapest, P.O.B. 127 H--1364, Hungary, e-mail:
major\@renyi.hu}
\medskip
 
{\narrower \noindent {\it Summary:}\/ In this paper a result of
Lata{\l}a about the tail behaviour of Gaussian polynomials will be
discussed. Lata{\l}a proved an interesting result about this problem
in paper~[2]. But his proof applied an incorrect statement at a
crucial point. Hence the question may arise whether the main result
of paper~[2] is valid. The goal of this paper is to settle this
problem by presenting such a proof where the application of the
erroneous statement is avoided. I discuss the proofs in detail even
at the price of a longer text and try to give such an explanation
that reveals the ideas behind them better than the original paper.
\par}

 
\beginsection 1. Introduction. Formulation of the main results.
 
In this paper the following problem studied in Lata{\l}a's paper~[2]
will be revisited.
 
Let us have a multilinear form
$$
\aligned
A(u_1,\dots,u_d)&=A(d)(u_1,\dots,u_d)\\
&=\sum_{(i_1,\dots,i_d)\colon\;1\le i_j\le n_j,\;1\le j\le d}
a(i_1,\dots,i_d) u_1(i_1)\cdots u_d(i_d)
\endaligned \tag1.1
$$
of order~$d$ in the space of vectors $(u_1,\dots,u_d)$ where
$u_j=(u_j(1),\dots,u_j(n_j))\in R^{n_j}$, and $R^{n_j}$ is the
Euclidean space with some prescribed dimension $n_j$, $1\le j\le d$.
The set of real numbers $A(d)=A(d|n_1,\dots,n_d)=
\{a(i_1,\dots,i_d),\; 1\le i_j\le n_j,\; 1\le j\le d\}$ are also
prescribed in this formula.
 
Beside this, let us also have $d$ independent standard Gaussian
random vectors $G_j=(g_j(1),\dots,g_j(n_j))$ of dimension $n_j$,
$1\le j\le d$, and define with the help of the multilinear form~(1.1)
and these Gaussian random vectors the Gaussian random polynomial
$$
Y(A)=Y(A(d))=\sum_{(i_1,\dots,i_d)\colon\;1\le i_j\le n_j,\;1\le j\le d}
a(i_1,\dots,i_d) g_1(i_1)\dots g_d(i_d) \tag1.2
$$
of order $d$. We want to give a good estimate on the tail
distribution $P(|Y(A)|>x)$ for all $x>0$ under appropriate
conditions on the multilinear form $A(d)(\cdot)$ defined in~(1.1).
Naturally, it belongs to the problem to find the right conditions
under which useful results can be proved.
 
Some estimates can be proved about the tail distribution of
Gaussian polynomials and so-called degenerate $U$-statistics under
the condition that their variance is bounded by a known constant,
(see~[5]), and these results are in a sense sharp. On the other hand,
they can be improved if we have some additional useful
information about the behaviour of the multi-linear  form~(1.1).
Lata{\l}a proved an interesting result in this direction in paper~[2].
He found the right conditions under which a good estimate can be given
about the tail-distribution $P(|Y(A)|>x)$. Similar questions can be
also asked about degenerate $U$-statistics, and Adamczak proved in~[1]
some results in this direction. But the essential step in the study of
such problems is to find the proof (and formulation) of the right
estimates for the tail distribution of Gaussian polynomials. The
adaptation of such results to $U$-statistics is rather a technical
problem.
 
Hence I restrict my attention to Lata{\l}a's work. I discuss its
proof and present a version of it, because I found an error in
paper~[2] that caused serious problems for me. For a long time I
have even doubted the validity of the main result in~[2]. My problems
were related to the proof of Theorem~3 in~[2]. It was based
on a backward induction procedure with respect to a parameter~$l$.
The induction steps when we turn from $l+1$ to $l$ were explained
for all parameters $l\ge1$. But the final step when we turn from 
$l=1$ to $l=0$ was not considered in the proof. Moreover,
the arguments of the paper do not work in this case, and as a
consequence the proof of Theorem~3 is invalid. At the end of
Section~8 I discuss this problem in a remark in more detail.
 
The above mentioned error seems to be crucial. I believe that not
only the proof but even the formulation of Theorem~3 is erroneous.
Since the proof of the main result of paper~[2] depends heavily
on this theorem the question arises whether this result holds. It
demanded  much work from me to clarify this problem. Finally I
found a correct proof of the main result of paper~[2] which does not
apply Theorem~3 of~[2]. I present it in this paper. Beside this I
also try to explain its main ideas.
 
To formulate Lata{\l}a's result I introduce some notations.
Let us define the linear functional
$$
A(v)=A(d,v)=\sum_{(i_1,\dots,i_d)\colon\;1\le i_j\le n_j,\;1\le j\le d}
a(i_1,\dots,i_d) v(i_1,\cdots, i_d) \tag1.3
$$
in the space of all real valued functions $v(i_1,\dots,i_d)$
defined on the $n$-tuples $(i_1,\dots,i_d)$, $1\le i_j\le n_j$,
$1\le j\le d$, where the coefficients $a(i_1,\dots,i_d)$
agree with those given in~(1.1)
 
Let us also introduce the class $\Cal P=\Cal P_d$ consisting of all
partitions of the set $\{1,\dots,d\}$. We shall define a class of
finite series of functions with the help of these partitions, and
the conditions of Lata{\l}a's result will be formulated with their
help. To avoid some repetitions in further discussions I define
these quantities in a slightly more general form.
 
Let us have a finite subset $K\subset\{1,2,\dots,\}$ of the positive
integers together with a function $b_K(i_j,\,j\in K)$, 
$1\le i_j\le n_j, \,j\in K$, and the numbers $n_j$, $j\in K$, which
tell what values the arguments of the function $B_K(i_j,\,j\in K$,
$1\le j\le n_j,\, j\in K)$ can take. We  define with their help, 
similarly to the quantity $A(v)$, the linear functional
$$
B_K(v)=\sum_{(i_j,\,j\in K)\colon\;1\le i_j\le n_j,\;j\in K}
b_K(i_j,\;1\le i_j\le n_j,\, j\in K) v(i_j,\;1\le i_j\le n_j,\, j\in K)
\tag1.4
$$
on the space of functions $v(i_j,\;1\le i_j\le n_j,\,j\in K)$.

Let $\Cal P(K)$ denote the set of all partitions of the set $K$, and
given a partition $P=\{A_1,\dots,A_s\}\in\Cal P(K)$ of $s$ elements
together with the positive integers $n_j$, $j\in K$, appearing in
the definition of the sets $P(K)$ let us define with their help the
following set $\Cal G_P$ of sequences of functions $(v_1,v_2,\dots,v_s)$:
$$
\aligned
\Cal G_P&=\biggl\{
(v_1(i_j,\; 1\le i_j\le n_j, j\in A_1),
\dots, v_s(i_j,\; 1\le i_j\le n_j, j\in A_s))\colon \\
&\qquad\qquad \sum_{(i_1,\dots,i_j)\colon\; 1\le i_j\le n_j,\, j\in A_r}
v_r^2(i_j,\,j\in A_r)\le1 \quad \text{for all } 1\le r\le s \biggr\}
\endaligned \tag1.5
$$
if $P=\{A_1,\dots,A_s\}\in\Cal P(K)$. Let us have a linear functional
$B_K(v)$ of the form~(1.4) together with the coefficients $b_K(\cdot)$
taking part in its definition. Then we define with the help of the
class of functions $\Cal G_P$ defined in (1.5) the following quantity
$V(P,B_K)$ for all partitions $P\in\Cal P(K)$.
$$
\aligned
&V(P,B_K)=V(P,b_K(\cdot)) \\
&\qquad =\sup_{(v_1,\dots,v_s)\in\Cal G_P}
\sum b_K(i_j,\;1\le i_j\le n_j,\,j\in K)
\prod_{1\le r\le s}v_r(i_j,\,1\le i_j\le n_j,\, j\in A_r).
\endaligned \tag1.6
$$
for a partition $P=\{A_1,\dots,A_s\}\in\Cal P(K)$. In this formula
the same coefficients $b_K(i_j,\;1\le i_j\le n_j,\,j\in K)$ appear
as in~(1.4).
 
Given a partition $P=\{A_1,\dots,A_s\}\in\Cal P(K)$ let $|P|=s$ denote
its cardinality. In the remaining part of this section I restrict
my attention to partitions $P\in\Cal P_d$ of the set $\{1,\dots,d\}$
and to the case when the linear functional $A(v)$ defined in
formula~(1.3) is considered. In this case the quantity introduced
in~(1.6) will be denoted as $V(P,A)=V(P,(a(\cdot))$. Let us define
with its help the numbers
$$
\alpha_s=\alpha_s(A)
=\sup_{P\colon\;P\in\Cal P_d,\,|P|=s}V(P,A)
\quad\text{for all } 1\le s\le d. \tag1.7
$$
The main result of Lata{\l}a we discuss in this paper can be formulated
with the help of the quantities $\alpha_s$, $1\le s\le d$, introduced
in~(1.7). It states the following inequalities.
 
\medskip\noindent
{\bf Theorem~1.} {\it The moments of the Gaussian random polynomial
$Y(A(d))$ defined in formula~(1.2) satisfy the inequality
$$
E(Y(A(d)))^{2M}\le \(C(d)\max_{1\le s\le d}(M^{s/2}\alpha_s)\)^{2M} \tag1.8
$$
for all $d\ge2$ and $M=1,2,\dots$ with the quantities $\alpha_s$
defined in~(1.7) and a constant $C(d)$ depending only on the
order~$d$ of the Gaussian polynomial~$Y(A(d))$. As a consequence,
$$
P(|Y(A(d))|>x)\le C(d)\exp\left\{-\frac1{C(d)}
\min_{1\le s\le d}\(\frac x{\alpha_s}\)^{2/s}\right\} \tag1.9
$$
for all $d\ge2$ and $x>0$ with some constant $C(d)$ depending
only on~$d$.}
 
\medskip\noindent
{\it Remark 1.}\/ Lata{\l}a's paper also contains a similar lower
bound for the moments and probabilities in~(1.8) and~(1.9). These
bounds state that the estimates in this formulas are essentially sharp,
only the value of the parameter~$C(d)$ can be improved in them. The
proof of these lower bounds is considerably simpler. Since their
proof in~[2] is correct, I shall omit their discussion.
 
\medskip\noindent
{\it Remark 2.}\/ In the subsequent estimations some constants $C$,
$C_1$, $C(d)$ etc. will appear in different formulas. The same letter
may denote different constants in different formulas. It will be
important that these constants are universal, depending at least of
the order $d$ of the Gaussian polynomial we are considering. There
will be some places in our discussion where the relation between
constants in different formulas have to be investigated. The necessary
considerations will be taken at these points.
 
\medskip\noindent
{\it Remark 3.}\/ The dimension $n_j$ of the Euclidean spaces $R^{n_j}$
where the appropriate vectors take their values plays no role in
our considerations. It is exploited in some arguments that they are
finite, but their value will be not important for us. At several
points where it makes no problem I shall omit the parameters~$n_j$
from the formulas. By means of some limiting procedure one can get
results in infinite dimensional spaces, but this will be not done here.
 
\medskip
I formulate a formally weaker version of Theorem~1 in the following
Theorem~1A. But actually, as I shall show these two results are
equivalent. Since Theorem~1A is technically simpler, this result
will be proved.
 
\medskip\noindent
{\bf Theorem~1A.} {\it Let the Gaussian polynomial $Y(A(d))$, $d\ge2$,
defined in~(1.2) be such that the expressions $\alpha_s$,
$1\le s\le d$, defined in~(1.7) satisfy the inequality
$$
\alpha_s=\alpha_s(A)\le M^{-(s-1)/2} \quad\text{for all } 1\le s\le d \tag1.10
$$
with some positive integer $M$. Then
$$
EY(A(d))^{2M}\le C(d)^MM^M \tag1.11
$$
with a constant $C(d)>0$ depending only on the order~$d$ of the
Gaussian polynomial~$Y(A(d))$.}
 
\medskip
Theorem~1A states that if a Gaussian polynomial $Y(A(d))$ satisfies
condition~(1.10) then its $2M$-th moment satisfies such an estimate
as the $2M$-th moment of a standard normal random variables multiplied
by a constant.
 
\medskip\noindent
{\it The deduction of Theorem~1 from Theorem~1A.}\/ Let us consider
the random variable $Y(A(d))$ and the number~$2M$ which is the moment
we consider in formula~(1.8). Let us define with their help the constant
$D(M)=\max\limits_{1\le s\le d}(M^{(s-1)/2}\alpha_s)$ and
introduce the Gaussian polynomial $D(M)^{-1}Y(A(d))$ defined in
formula~(1.2) with coefficients $D(M)^{-1}a(i_1,\dots,i_d)$. This
polynomial satisfies relation~(1.10), hence by Theorem~1A relation~(1.11)
also holds for it. This means that $EY(A(d))^{2M}\le (C(d)D(M)^2M)^M$
which is equivalent to relation~(1.8) in Theorem~1.
 
Relation~(1.9) follows from relation~(1.8) in the standard way. By
the Markov inequality $P(|Y(A(d))|\ge x)\le x^{-2M}EY(A(d))^{2M}$
for arbitrary $M=1,2,\dots$. Choose $M=\[\min\limits_{1\le s\le d}
\frac1{K C(d)}\frac x{\alpha_s}\]^{2/s}$
if $x\ge K C(d)\min\limits_{1\le s\le d}\alpha_s$, where 
$[\cdot]$ denotes integer part, $C(d)$ is the same constant which 
appears in (1.8), and $K=K(d)$ is a sufficiently large constant 
depending only on~$d$. In this case
we get from relation~(1.8) that $P(|Y(A(d))|\ge x)\le e^{-M}$ which
implies relation~(1.9) with the constant $K^2C(d)^2$ if 
$x\ge KC(d)\min\limits_{1\le s\le d}\alpha_s$. On the other hand, 
if $x\le KC(d)\min\limits_{1\le s\le d}\alpha_s$, and the constant
$K$ was chosen sufficiently large, then the right-hand side
of relation~(1.9) (with the previously chosen constant $K^2C^2(d)$ 
as the number `$C(d)$' in~(1.9)) is larger than~1. Hence 
relation~(1.9) holds also in this case.
 
\medskip
This paper consists of eight sections and an Appendix. In Section~2
the proof of Theorem~1A is reduced to a result called the Basic
estimate by means of a conditioning argument. In Section~3 this
Basic estimate is proved in the special case~$d=2$. In Section~4
a result of paper~[2] is recalled about the estimation of the
cardinality of an appropriate $\e$-net in a metric space with some
nice properties. In Section~5 a result called the Main inequality
is presented, and it is shown that the Basic estimate follows
from it. In Section~6 two results, Lemma~6.1 and Lemma~6.2 are
formulated. They provide a good partition of certain sets of
functions which play crucial role in the proof of the Main
inequality. The proof of these lemmas is based on some estimates
formulated in Lemma~6.3. Lemma~6.3 together with its proof is
also given in Section~6. Lemmas~6.1 and~6.2 are proved in
Section~7. Finally the Main inequality is proved in Section~8 by
means of the results in Section~6. Since in Section~4 I apply a
terminology essentially different from that  of~[2] I found better
not to refer to the original proofs of the results presented
here, but to describe them instead. This is done in the
Appendix. In such a way I wanted to make this paper self-contained.
 
The proofs of this paper apply several ideas of Paper~[2]. But
since the notation and the formulation of the results in these
two works are very different, and the main ideas in~[2] are
presented in a rather hidden way I only explain which results
of these two paper correspond to each other.
 
\beginsection 2. The application of a conditioning argument.
 
In this section a conditioning argument is applied to reduce the
proof of Theorem~1A to the verification of a result called the
Basic estimate.
 
To carry out this conditioning argument let us define the
Gaussian random vector
$$
Y_d(u)=Y_d(u,A)=\sum_{(i_1,\dots,i_d)\colon\;
1\le i_j\le n_j,\;1\le j\le d}a(i_1,\dots,i_d)
u_1(i_1)\dots u_{d-1}(i_{d-1}) g_d(i_d) \tag2.1
$$
for all vectors $u=(u_1,\dots,u_{d-1})$, $u_j=(u_j(1),\dots,u_j(n_j))$,
$1\le j\le d-1$, and a standard Gaussian vector
$G_d=(g_d(i_1),\dots,g_d(n_d))$. The coefficients $a(i_1,\dots,i_d)$
in formulas (1.1) and (2.1) are the same. Actually in formula (2.1)
we took the multilinear form~(1.1) and replaced the vector $u_d$ by
the standard normal random vector~$G_d$ in it.
 
We want to estimate the moments of the random variables $Y(A(d))$
introduced in~(1.2). This can be done by means of the following
conditioning argument.
$$
\align
&E(Y(A(d))^{2M}|g_d(1)=u_d(1),\dots,g_d(n_d)=u_d(n_d))\\
&\qquad=E\(\sum_{(i_1,\dots,i_d)\colon\;
1\le i_j\le n_j,\;1\le j\le d}a(i_1,\dots,i_d)
g_1(i_1)\dots g_{d-1}(i_{d-1})u_d(i_d)\)^{2M},
\endalign
$$
hence
$$
EY(A(d))^{2M}=EY(A(d),M,G_d) \tag2.2
$$
with
$$
\aligned
&Y(A(d),M,u_d)\\
&\quad =E\[\sum_{i_d=1}^{n_d}
\(  \! \sum_{(i_1,\dots,i_{d-1})\colon\;1\le i_j\le n_j,\, 1\le j\le d-1}
\!\!\! a(i_1,\dots,i_d)g_1(i_1)\dots g_{d-1}(i_{d-1})\)u_d(i_d)\]^{2M},\!\!\!\\
\endaligned
$$
or in an equivalent form
$$
\aligned
&Y(A(d),M,u_d)\\
&\qquad =E\(\sum_{(i_1,\dots,i_{d-1})
\colon\;1\le i_j\le n_j,\, 1\le j\le d-1}
b_{u_d}(i_1,\dots,i_{d-1}) g_1(i_1)\dots g_{d-1}(i_{d-1})\)^{2M},
\endaligned \tag2.3
$$
with
$$
b_{u_d}(i_1,\dots,i_{d-1})=\sum_{i_d=1}^{n_d} a(i_1,\dots,i_d)u_d(i_d),
\tag2.4
$$
where $u_d=(u_d(1),\dots,u_d(n_d))$ is an arbitrary vector in $R^{n_d}$.
 
Next I formulate a result called the Basic estimate. Its proof will
be the main subject of the subsequent sections. Here I prove that
Theorem~1A follows from it. To formulate the Basic estimate first I
introduce the following quantity.
$$
Z_d=Z_d(A)=\sup_{u=(u_1,\dots,u_{d-1})\colon\; u_j\in B^{n_j},
\,1\le j\le d-1} Y_d(u), \tag2.5
$$
where the (Gaussian) random variables  $Y_d(u)$ were defined in~(2.1).
Here and in the subsequent part of the paper $B^{n}$ denotes the unit
ball in the Euclidean space $R^n$ with the usual Euclidean norm, i.e.
$B^n=\{(u(1),\dots,u(n))\colon\;\summ_{j=1}^n u(j)^2\le1\}$. It will
be shown with the help of the previous calculations that Theorem~1A
follows from the following result.
 
\medskip\noindent
{\bf Basic estimate.} {\it If the linear form $A(v)$, $d\ge2$,
introduced in~(1.3) is such that the quantities $\alpha_s$
defined in~(1.7) satisfy the condition~(1.10) with some positive
integer $M$, i.e. $\alpha_s=\alpha_s(A)\le M^{-(s-1)/2}$ for all
$1\le s\le d$, then the estimate
$$
EZ_d^{2M}=EZ_d(A)^{2M}\le C^M M^{-(d-2)M} \tag2.6
$$
holds with a constant $C=C(d)$ depending only on $d$.}
 
\medskip\noindent
{\it Remark.} The above formulated Basic estimate is closely
related to Theorem~2 in~[2]. The main difference between them is
that Theorem~2 in~[2] gives an estimate only for the expected
value $EZ_d(A)$ of $Z_d(A)$ and not for its higher moments. Thus
our result is, --- at least formally, --- sharper. But actually
estimate~(2.6) follows from the result of~[2] and an important
concentration inequality of Ledoux about the supremum of Gaussian
random variables. This result will be recalled in Section~3. The
reason for the present formulation of the Basic estimate was that
I wanted to show that the so-called chaining argument applied in
its proof also supplies the estimate~(2.6) for $d\ge3$, i.e. we
do not need Ledoux's inequality in this case. Surprisingly, we
need it just in the simplest case $d=2$, when the proof is given
by means of a simple and natural direct calculation instead of
the chaining argument.

\medskip
We shall estimate $EY(A(d))^{2M}$ with the help of relations~(2.2)
and~(2.3) by induction with respect to $d$ for all $d\ge2$.
Let us first consider the case $d=2$.
 
If the linear form $A(2)(u_1,u_2)$ in~(1.1) (with $d=2$) is defined
with the help of a set of numbers
$\{a(i,j)\; 1\le i\le n_1,\, 1\le j\le n_2\}$, then we can write
$$ \allowdisplaybreaks
\align
&Y(A(2),M,u_2)=
E\[\sum_{i=1}^{n_1}\(\sum_{j=1}^{n_2} a(i,j) u_2(j)\)g_1(i)\]^{2M}\\
&\qquad =1\cdot3\cdot\cdots\cdot (2M-1)\(E\[\sum_{i=1}^{n_1}
\(\sum_{j=1}^{n_2} a(i,j)u_2(j)\) g_1(i)\]^2\)^{M} \\
&\qquad =1\cdot3\cdot\cdots\cdot (2M-1)\(\sum_{i=1}^{n_1}
\(\sum_{j=1}^{n_2} a(i,j)u_2(j)\)^2\)^{M}  \tag2.7 \\
&\qquad=1\cdot3\cdot\cdots\cdot (2M-1)
\(\sup_{u_1=(u_1(1),\dots,u_1(n_1))\colon\; u_1\in B^{n_1}}
\sum_{i=1}^{n_1}\sum_{j=1}^{n_2} a(i,j) u_1(i) u_2(j)\)^{2M},
\endalign
$$
where $u_2=(u_2(1),\dots,u_2(n_2)\in R^{n_2}$,
$u_1=(u_1(1),\dots,u_1(n_1))\in B^{n_1}$, and $B^{n_1}$
denotes the unit ball of the Euclidean space $R^{n_1}$, i.e.
we demand that $\sum\limits_{i=1}^{n_1}u_1(i)^2\le1$. By
relations~(2.2), (2.7), the definition of the quantity $Z_d(A)$ 
and the Basic estimate
$$
EY_2(A(2))^{2M}\le (2M)^MEZ_2(A(2))^{2M}\le CM^M
$$
if $\alpha_1(A)\le1$ and $\alpha_2(A)\le M^{-1/2}$, i.e. if the
conditions of the Basic estimate hold for $d=2$. Thus we have
proved Theorem~1A with the help of the Basic estimate in the
case~$d=2$.
 
In the case $d\ge3$ Theorem~1A will be proved by means of
induction. During this induction procedure we assume that
Theorem~1A holds for $2\le d'\le d-1$, and the Basic estimate
holds for $2\le d'\le d$.
 
First the expression $Y(A(d),M,u_d)$  will be estimated. This
expression, defined in (2.3) is the $2M$-th moment of a Gaussian
polynomial of order $d-1$. It is defined similarly to $Y_d(u)$
introduced in formula~(2.1) only with the coefficients
$b_{u_d}(i_1,\dots,i_{d-1})$ introduced in~(2.4) instead of
$a(i_1,\dots,i_d)$. Hence, as we shall show, they satisfy the
following inequality.
$$
\aligned
Y(A(d),M,u_d)&\le \max_{P\in\Cal P_{d-1}}\(V(P,B_{u_d})^2M^{(|P|-1)}\)^M
(CM)^M\\
&\le C^M \sum_{P\in\Cal P_{d-1}} V(P,B_{u_d})^{2M} M^{|P|M}
\endaligned \tag2.8
$$
with some constant $C=C(d)$, where $V(P,B_{u_d})$ was defined
in~(1.6) for partitions $P\in\Cal P_{d-1}$ i.e. $K=\{1,\dots,d-1\}$,
and the numbers $b_{u_d}(i_1,\dots,i_{d-1})$ introduced in~(2.4) play
the role of the coefficients $b_K(\cdot)$ in formulas~(1.4) and~(1.6).
 
Indeed, the expression $\frac{Y(A(d),M,u_d)}
{\max\limits_{P\in\Cal P_{d-1}}(V(P,B_{u_d})M^{(|P|-1)/2})^{2M}}$
equals the $2M$-th moment of such a Gaussian polynomial which satisfies
the conditions of Theorem~1A with parameter $d-1$. Hence Theorem~1A
with parameter $d-1$ (which holds by our induction hypothesis)
implies the first inequality in~(2.8). The second inequality of~(2.8)
is obvious.

By relations (2.2) and (2.8)
$$
EY(A(d))^{2M}\le C^M\sum_{P\in\Cal P_{d-1}}EV(P,B_{G_d})^{2M}M^{|P|M},
$$
where $V(P,B_{G_d})$ is the random variable we get by replacing the
vector $u_d$ by the random vector $G_d=(g_d(1),\dots,g_d(n_d))$ in
the expression $V(P,B_{u_d})$. Hence to complete the proof of the
Theorem~1A it is enough to show that under the conditions of
Theorem~1A
$$
EV(P,B_{G_d})^{2M}\le C^M M^{-(|P|-1)M} \quad \text{for all }
P\in\Cal P(\{1,\dots,d-1\}) \tag2.9
$$
with a constant $C=C(d)$. This result can be proved with the help
of the Basic estimate.

To prove formula~(2.9) take a partition
$P=\{A_1,\dots,A_s\}\in\Cal P_{d-1}$ with $|P|=s$ elements.
With such a choice
$$
V(P,B_{G_d})=\sup_{(v_1,\dots,v_s)\in\Cal G_P}
\sum_{(i_1,\dots,i_d)} a(i_1,\dots,i_d)
\prod_{r=1}^s v_r(i_j,\;j\in A_r) g_d(i_d). \tag2.10
$$
In formula~(2.10) the class of functions $\Cal G_P$ where the
supremum is taken is defined in~(1.5) with the partition~$P$ we
have fixed, and $(g_d(1)$,\dots, $g_d(n_d))$ is an $n_d$
dimensional standard normal vector. The $2M$-th moment of the 
right-hand side expression in~(2.10) can be bounded by means of 
the Basic estimate with $s+1=|P|+1\le d$ parameters (i.e. the 
number $|P|+1$ takes the role of the parameter $d$ in this case) 
if the vectors $(i_j,\,j\in A_r)$, $A_r\in P$, are considered as 
one variable for all $1\le r\le s$. The condition of the Basic
estimate formulated in~(1.10) holds with such a choice, and we 
get inequality~(2.9) in such a way.

We have reduced the problem we want to solve to the proof of an inequality
formulated in the Basic estimate, where certain moments of a supremum
$\sup\limits_{u\in B^{n_1}\times\dots\times B^{n_{d-1}}}Y_d(u)$ of
Gaussian random variables are bounded. The random variables $Y_d(u)$
in this formula were defined in~(2.1), and $B^n$ denotes the unit
ball in $R^n$. In the study of such problems it is worth introducing
the metric $\rho(u,v)=[E(Y_d(u)-Y_d(v))^2]^{1/2}$ on the parameter
set of the random variables we are considering. This led to the
definition of the following pseudometric $\rho_\alpha$ in the
space $R^{n_1}\times\cdots\times R^{n_{d-1}}$.
$$
\aligned
\rho_\alpha(u,v)&=\rho_\alpha((u_1,\dots,u_{d-1}),(v_1,\dots,v_{d-1}))\\
&=[E(Y_d(u)-Y_d(v))^2]^{1/2}
=E\bigg(E\biggl[\sum_{1\le i_j\le n_j,\,1\le j\le d}
a(i_1,\dots,i_d)\\
&\qquad\qquad\qquad (u_1(i_1)\cdots u_{d-1}(i_{d-1})
-v_1(i_1)\cdots v_{d-1}(i_{d-1}))g(i_d)\biggr]^2\biggr)^{1/2} \\
&=\biggl(\sum_{1\le i_d\le n_d}
\biggr[\sum_{1\le i_j\le n_j,\,1\le j\le d-1} a(i_1,\dots,i_d) \\
&\qquad\qquad\qquad (u_1(i_1)\cdots u_{d-1}(i_{d-1})
-v_1(i_1)\cdots v_{d-1}(i_{d-1}))\biggr]^2\biggr)^{1/2}
\endaligned \tag2.11
$$
for all pairs of vectors $u=(u_1,\dots,u_{d-1})$ and
$v=(v_1,\dots,v_{d-1})$, $u_j\in R^{n_j}$, $v_j\in R^{n_j}$,
$1\le j\le d-1$.
 
It is useful to give a different characterization of the
above introduce metric $\rho_\alpha$. For this goal let us
define the pseudonorm $\alpha$
$$
\aligned
\alpha(v)&=
\alpha_d(v)=\alpha_d(v(i_1,\dots,i_{d-1}))\\
&=\biggl[\sum_{1\le i_d\le n_d}
\biggl(\sum_{1\le i_j\le n_j,\,1\le j\le d-1} a(i_1,\dots,i_d)
v(i_1,\cdots,i_{d-1})\biggr)^2\biggr]^{1/2}
\endaligned \tag2.12
$$
in the linear space of the functions
$v=v(i_1,\dots,i_{d-1})$, $1\le i_j\le n_j$, $1\le j\le d-1$.
Clearly,
$$
\rho_\alpha((u_1,\dots,u_{d-1}),(v_1,\dots,v_{d-1}))
=\alpha_d(u_1\otimes\cdots\otimes u_{d-1}-v_1\otimes\cdots\otimes v_{d-1})
\tag2.13
$$
where the function $u_1\otimes\cdots\otimes u_{d-1}$ with
arguments $(i_1,\dots,i_{d-1})$, $1\le i_j\le n_j$ for all
$1\le j\le d-1$ is defined as
$u_1\otimes\cdots\otimes u_{d-1}(i_1,\dots,i_{d-1})=u_1(i_1)\cdots
u_{d-1}(i_{d-1})$, and $v_1\otimes\cdots\otimes v_{d-1}$ is defined
similarly.
 
The above representation of the metric $\rho_\alpha$ turned
out to be useful. In the study of the Basic estimate we have
to find a good $\e$-net for certain subsets of
$B^{n_1}\times\cdots\times B^{n_{d-1}}$ with respect to the
metric $\rho_\alpha$ for small $\e>0$. The representation of
the metric $\rho_\alpha$ by formulas~(2.12) and~(2.13) may help
in finding good $\e$-nets. This question will be discussed in
detail in the subsequent sections. But before doing it I prove
the Basic estimate together with some related results we need
in our discussion in the special case $d=2$. This case is
considered separately, because the formulation of the results
and their proof for $d=2$ are slightly different from those in
the general case.
 
\beginsection 3. The proof for Gaussian polynomials of order 2.
 
In this section the Basic estimate will be proved for Gaussian
polynomials of order $d=2$. It will be proved as the consequence
of a more general result called the Main inequality in the case
$d=2$. A result called the Main inequality will be formulated in
Section~5 for all dimensions $d\ge3$. The crucial point in the
proof of Theorem~1A is the verification of this result. The Main
inequality in the case $d=2$ formulated in this section can be
considered as a version of this result. But there are some
differences between their formulation, and they must be
considered separately. The Basic estimate for $d=2$ could have
been proved directly. I prove it with the help of the Main
inequality in the case $d=2$, because the latter result is also
needed in the discussion of the case $d\ge3$. To formulate it
I introduce some notations.
 
We shall work with some expressions $A(v)$ and $Y_2$ which are
the quantities defined in~(1.3) and~(2.1) in the special case $d=2$.
Let us write them down in more detail.
 
These terms depend on a set of numbers
$A=A(2)=\{a(i,j),\;1\le i\le n_1,\,1\le j\le n_2\}$. The first of
them is the linear functional
$$
A(v)=A(2,v)=\summ_{i,j} a(i,j)v(i,j)
$$
in the space of all functions $v(i,j)$ with arguments
$1\le i\le n_1$, $1\le j\le n_2$. This is the expression~(1.3)
in the case $d=2$. The expression~(2.1) can be written as
$$
Y(u)=Y_2(u)=\sum_{i,j}a(i,j)u(i)g_2(j),
$$
with $u=(u(1),\dots,u(n_1))$, where $(g_2(1),\dots,g_2(n_2))$ is 
a standard normal random vector.
 
Let us observe that in the case $d=2$ the quantity $\alpha_1(A)$
defined in~(1.7) can be written as
$$
\alpha_1(A)=\sup_{v(i,j)\colon\;\summ_{i,j}v(i,j)^2\le1}
\sum_{i,j} a(i,j)v(i,j)
=\(\sum_{i,j} a(i,j)^2\)^{1/2}. \tag3.1
$$
Let us also introduce the function
$$
\alpha_2(u)=\[\sum_j\(\sum_{i} a(i,j)u(i)\)^2\]^{1/2}
=\sup_{v=(v(1),\dots,v(n_2))\colon\; \summ_j v(j)^2\le1}
\sum_{i,j} a(i,j)u(i)v(j)
$$
for all vectors $u=(u(1),\dots,u(n_1))\in R^{n_1}$.
 
Let us fix some positive integer $M$, and define for all $N\ge0$
the following subset $U_N=U_N(M)$ of $R^{n_1}$.
$$
U_N=U_N(M)=\{u=(u(1),\dots,u(n_1))\colon\;  u\in B^{n_1},
\text{ and } \alpha_2(u)\le 2^{-N}M^{-1/2}\}. \tag3.2
$$
I formulate with the help of the above notations the following
result.
 
\medskip\noindent
{\bf  The Main inequality in the case $d=2$.} {\it Let
$\alpha_1(A)\le1$. Then the inequality
$$
E\[\sup_{u\colon\; u\in U_N}Y(u)\]^{2^{2(N+A)}M}
\le (C\cdot2^A)^{2^{2(N+A)}M} \tag3.3
$$
holds with the sets $U_N$ defined in~(3.2) for all integers 
$N\ge0$, $M\ge1$ and $A\ge 1$ with $C=2$.}
 
\medskip\noindent
{\it Proof of the Main inequality in the case $d=2$.}\/
This result will be proved with the help of the concentration
inequality of Ledoux about the supremum of Gaussian random
variables. (See~[3] Theorem~7.1.) First I show that under the
condition $\alpha_1(A)\le1$
$$
E\(\sup_{u=(u(1),\dots,u(n_1))\colon\; \summ_{i=1}^{n_1}u(i)^2\le1}Y(u)\)
\le 1. \tag3.4
$$
Indeed, for all $\oo\in\Omega$
$$
\sup_{ \summ_i u(i)^2\le1}
\sum_{i,j}a(i,j)u(i)g(j)(\oo)=\[\sum_i\(\sum_j a(i,j)g(j)(\oo)\)^2\]^{1/2},
$$
since the above expression takes its supremum at the value
$$
u(i)=\frac{\summ_j a(i,j)g(j)(\oo)}
{\[\summ_i\(\summ_j a(i,j)g(j)(\oo)\)^2\]^{1/2}}, \qquad 1\le i\le n_1.
$$
Hence by the Schwarz inequality and relation~(3.1)
$$
\align
&E\(\sup_{u=(u(1),\dots,u(n_1))\colon\; \summ_i u(i)^2\le1} Y(u)\)
=E\(\sup_{\summ_i u(i) ^2\le1}
\sum_{i,j}a(i,j)u(i)g(j)\)\\
&\qquad=E\[\sum_i\(\sum_j a(i,j)g(j)\)^2\]^{1/2}
\le \[E\sum_i\(\sum_j a(i,j)g(j)\)^2\]^{1/2}\\
&\qquad=\(\sum_{i,j}a(i,j)^2\)^{1/2}=\alpha_1(A)\le1.
\endalign
$$
On the other hand $EY(u)=0$ and $EY(u)^2=\alpha_2(u)^2\le 2^{-2N}M^{-1}$,
for all $u\in U_N$. Hence Ledoux's concentration inequality
(see formula ~7.4 in~[3]) implies that
$$
P\(\sup_{u\in U_N} \left|Y(u)-E\sup_{u\in U_N}Y(u)\right|\ge x\)
\le 2e^{-2^{-2N-1}Mx^2} \quad\text{for all } x\ge0.
$$
The above inequality with partial integration yield for all $R\ge2$ that
$$
\align
&E\sup_{u\in U_N} \left|Y(u)-E\sup_{u\in U_N}Y(u)\right|^{2R}
\le \int_0^\infty 2e^{-2^{2N-1}Mx^2}\,dx^{2R}\\
&\quad =4R\cdot2^{-2NR}M^{-R}\int_0^\infty x^{2R-1}e^{-x^2/2}\,dx
=4R\cdot2^{-2NR}M^{-R}(2R-2)(2R-4)\cdots2\\
&\quad\le(2RM^{-1})^R2^{-2NR}=(2RM^{-1}2^{-2N})^R.
\endalign
$$
Relation~(3.3) follows from the above inequality with the choice
$2R=2^{2(N+A)}M$, $N\ge0$, $M\ge1$, $A\ge1$, and the inequality
$E\supp_{u\in U_N} Y(u)\le1$ which is a consequence of relation~(3.4).
 
\medskip\noindent
{\it Proof of the Basic estimate for $d=2$.}\/ Let us apply the Main
inequality in the case $d=2$ with $N=0$ and $A=1$. Since the
conditions of the Basic estimate for $d=2$ contain the inequality
$\alpha_2(A)=\sup\limits_{u\in B^{n_1}}\alpha_2(u)\le M^{-1/2}$
the set $U_0$ agrees with the unit ball $B^{n_1}$.
Hence the Schwarz inequality and relation~(3.3) with the choice
$N=0$ and $A=1$ yield the estimate
$$
\align
E\[\sup_{u\colon\; u\in B^{n_1}}\sum_{i,j} a(i,j)u(i)g(j)\]^{2M}
&\le \(E\[\sup_{u\colon\; u\in U_0}\sum_{i,j} a(i,j)u(i)g(j)\]
^{4M}\)^{1/2} \\
&\le 4^{4M/2}=2^{4M}.
\endalign
$$
The Basic estimate for $d=2$ (with $C=16$ in formula~(2.6)) is proved.
 
\beginsection 4. Estimates on the cardinality of $\e$-nets with
respect to nice metrics.

In the Basic estimate the moments of the supremum of a class of
Gaussian random variables are estimated. In such problems it is
worth introducing a natural metric on the set of parameters of
the random variables we are considering, by defining the
distance of two points in the parameter space as the square root
of the variance of the difference of the corresponding random
variables. It is also useful to find such a subset of the parameter
space with relatively small cardinality which is dense with
respect to this metric. Such an approach leads to the
formulation of the following problem.
 
Given a pseudometric space $(X,\rho)$ together with a subset
$X_0\subset X$ we want to find for all $\e>0$ an $\e$-net of
relatively small cardinality in the space $X_0$ with respect to
the metric $\rho$, i.e. we want to find a set
$\{x_1,\dots,x_N\}\subset X_0$ with a relatively small index $N$ for
which $\min\limits_{1\le j\le N}\rho(x_j,x)\le\e$ for all $x\in X_0$.
A good $\e$-net can be found by solving the following problem.
Let us define an appropriate probability measure $\mu$ in the space
$(X,\rho)$ and give a good lower bound on the probability
$\mu(\{y\colon\;y\in X,\;\rho(y,x)\le\e\})$ for all $x\in X_0$ and
$\e>0$.
 
Lata{\l}a presented two estimates of this kind in Lemmas~1 and~2
of his paper~[2]. In Lemma~1 that case is considered  when $X$ is
the $n$-dimensional Euclidean space $R^n$, $X_0$ is the unit ball
in this space with respect to the Euclidean metric, and the
pseudometric $\rho=\rho_\alpha$ is defined by means of a
pseudonorm $\alpha$ in $R^n$ in the usual way, i.e.\
$\rho_\alpha(x,y)=\alpha(x-y)$. Lemma~2 is a multi-linear version
of this result. Here the space $X$ is the product of some
Euclidean spaces. We embed it in the tensor product of these
Euclidean spaces in a natural way, and the metric $\rho_\alpha$
in $X$ is defined with the help of a pseudonorm in this tensor
product.
 
Since these results play an important role in our considerations
I recall them in this paper under the names Proposition~4.1 and
Proposition~4.2. I shall apply a notation different from~[2],
and it may be hard to compare the results formulated here with
their original version. Hence to make this paper self-contained
I present the proof of Lata{\l}a's results in an Appendix.
 
To formulate these results some notations have to be introduced.
We denote the unit ball in the $n$-dimensional Euclidean space
by~$B^n$. We introduce a probability measure $\mu_{n,t}$ depending
on a parameter~$t$ in the Euclidean space~$R^n$ in the following
way. Given some number $t>0$ let $\mu_{n,t}$ denote the
distribution of the random vector $tG=(tg_1,\dots,tg_n)$ in
$R^n$, where $g_1,\dots,g_n$ are independent standard normal
random variables.
 
\medskip\noindent
{\bf Proposition 4.1.}  {\it Let $\alpha_1$ and $\alpha_2$
be two pseudonorms in $R^n$, $t>0$ an arbitrary positive number,
$x\in B^n$ a vector in the unit ball of $R^n$ and
$G=(g_1,\dots,g_n)$ an $n$-dimensional standard normal vector.
Then
$$
\mu_{n,t}(\{y\colon\; y\in R^n,\;\alpha_1(y-x)\le 4E\alpha_1(tG),\;
\alpha_2(y-x)\le 4E\alpha_2(tG)\})\ge \frac12e^{-1/2t^2}
$$
with the above introduced probability measure $\mu_{n,t}$.}
 
\medskip\noindent
{\it Remark.}\/ In our applications it would be enough to
consider a simpler version of Proposition~4.1 where only one
pseudonorm $\alpha_1$ appears. We formulated a result with two
pseudonorm, because such a result is applied in the proof of
Proposition~4.2.
 
\medskip
To formulate Proposition~4.2 some additional notations have to be
introduced. Let us consider $d$ Euclidean spaces 
$R^{n_1},\dots,R^{n_d}$ of dimension $n_j$, $1\le j\le d$, their 
product $R^{n_1}\times\cdots\times R^{n_d}$ and their tensor 
product $R^{n_1}\otimes\cdots\otimes R^{n_d}$ with some pseudonorm 
$\alpha(\cdot)$ on the tensor product. We give an embedding of
the product $R^{n_1}\times\cdots\times R^{n_d}$  of these Euclidean
spaces into their tensor product and
define with its help a pseudometric $\rho_\alpha$ in the
product space $R^{n_1}\times\cdots\times R^{n_d}$
induced by the pseudonorm $\alpha$ on the tensor product
$R^{n_1}\otimes\cdots\otimes R^{n_d}$.
 
For the sake of simpler notations we shall represent the
Euclidean space $R^n$ as the space of the real valued functions
$x=(x(1),\dots,x(n))$ on the set $\{1,\dots,n\}$, the tensor
product $R^{n_1}\otimes\cdots\otimes R^{n_d}$ of
the Euclidean spaces $R^{n_j}$, $1\le j\le d$, as the space of
the real valued functions $v(i_1,\dots,i_d)$, defined on the
set of vectors  $(i_1,\dots,i_d)$, $1\le i_j\le n_j$,
$1\le j\le d$, and the product $R^{n_1}\times\cdots\times R^{n_d}$
as the space of all vectors $x=(x_1,\dots,x_d)$, whose elements
are real valued functions $x_j=(x_j(1),\dots,x_j(n_j))$ on
the sets $\{1,\dots,n_j\}$, $1\le j\le d$.
 
We embed the Euclidean space $R^{n_1}\times\cdots\times R^{n_d}$
in the tensor product $R^{n_1}\otimes\cdots\otimes R^{n_d}$
with the help of the map
$A(x)=A(x_1,\dots,x_d)=x_1\otimes\cdots\otimes x_d$ from the
Euclidean space $R^{n_1}\times\cdots\times R^{n_d}$ into the tensor
product $R^{n_1}\otimes\cdots\otimes R^{n_d}$,
where $x_1\otimes\cdots\otimes x_d$ is defined for a vector
$x=(x_1,\dots,x_d)\in R^{n_1}\times\cdots\times R^{n_d}$ by the
formula
$x_1\otimes\cdots\otimes x_d(i_1,\dots,i_d)=x_1(i_1)\cdots x_d(i_d)$
for all coordinates $(i_1,\dots,i_d)$ with $1\le i_j\le n_j$,
$1\le j\le d$.
 
Given a pseudonorm $\alpha$ on the tensor product
$R^{n_1}\otimes\cdots\otimes R^{n_d}$ define with
its help the pseudometric $\rho_\alpha$ in the space
$R^{n_1}\times\cdots\times R^{n_d}$ by the formula
$$
\rho_\alpha((x_1,\dots x_d),(y_1,\dots,y_d))
=\alpha(x_1\otimes\cdots\otimes x_d-y_1\otimes\cdots\otimes y_d) \tag4.1
$$
for all $x=(x_1,\dots,x_d)\in R^{n_1}\times\cdots\times R^{n_d}$ and
$y=(y_1,\dots,y_d)\in R^{n_1}\times\cdots\times R^{n_d}$. I shall
call this $\rho_\alpha$ the pseudometric induced by the
pseudonorm~$\alpha$.
 
Let us fix some $x=(x_1,\dots,x_d)\in B^{n_1}\times\cdots\times B^{n_d}$
in the product of the unit balls $B^{n_j}$ in $R^{n_j}$,
$1\le j\le d$. In Proposition~4.2 a good lower bound is given on
the probability of a small neighbourhood of such a point~$x$ with
respect to an appropriately defined probability measure. More
explicitly, the probability
$\mu_{n_1+\cdots+n_d,t}(y\colon\;y\in R^{n_1}\times\cdots\times R^{n_d},
\rho_\alpha(x,y)\le u)$ will be bounded from below for all numbers
$u>0$ with respect to an appropriately defined Gaussian measure
$\mu_{n_1+\cdots+n_d,t}$, where $\rho_\alpha$ is the pseudometric
in $R^{n_1}\times\cdots\times R^{n_d}$ induced by a
pseudonorm $\alpha$ in $R^{n_1}\otimes\cdots\otimes R^{n_d}$ in the
above way. To formulate this result some additional
notations will be introduced.
 
Let us consider $d$ independent standard normal vectors
$G_j=(g_j(1),\dots,g_j(n_j))$ of dimension $n_j$, $1\le j\le d$,
and for all $t>0$ let $\mu_{n_1+\cdots+n_d,t}$ denote the
distribution of the random vector $(tG_1,\dots,tG_d)$ in the
space $R^{n_1}\times\cdots\times R^{n_d}$. Given a
pseudonorm $\alpha$ on the tensor product
$R^{n_1}\otimes\cdots\otimes R^{n_d}$ of
the spaces $R^{n_j}$, $1\le j\le d$, a number $t>0$, some
set $I\subset\{1,\dots,d\}$, $I\neq\emptyset$ and a vector
$x=(x_1,\dots,x_d)\in R^{n_1}\times\cdots\times R^{n_d}$ we
define the quantity
$$
W_I^x(\alpha,t)=E\alpha(z_1\otimes\cdots\otimes z_d)\quad
\text{where $z_j=x_j$ if $j\notin I$ and $z_j=tG_j$ if $j\in I$} \tag4.2
$$
with the previously defined function $z_1\otimes\cdots\otimes z_d\in
R^{n_1}\otimes\cdots\otimes R^{n_d}$ for $(z_1,\dots,z_d)\in
R^{n_1}\times\cdots\times R^{n_d}$. In words, we take the function
$\alpha(x_1\otimes\dots\otimes x_d)$, replace the coordinates
$x_j\in R^{n_j}$ by $tG_j\in R^{n_j}$ for the indices
$j\in I$, and take the expected value of the random variable
obtained in such a way. With the help of the above quantities
we can formulate Proposition~4.2.
 
\medskip\noindent
{\bf Proposition~4.2.}  {\it Let us have a pseudometric
$\rho_\alpha$ in the product $R^{n_1}\times\cdots\times R^{n_d}$
of some Euclidean spaces $R^{n_j}$, $1\le j\le d$, induced by
a pseudonorm~$\alpha$ in their tensor product
$R^{n_1}\otimes\cdots\otimes R^{n_d}$. Fix some vector
$x=(x_1,\dots,x_d)\in B^{n_1}\times\cdots\times B^{n_d}$,
in the product of the unit balls $B^{n_j}$ in $R^{n_j}$,
$1\le j\le d$. The following inequality holds for such a
vector~$x$ and an arbitrary number $t>0$.
$$
\aligned
\mu_{n_1+\cdots+n_d,t}&\(\left\{y\colon\; y\in
R^{n_1}\times\cdots\times R^{n_d},\;
\rho_\alpha(x,y)\le \sum_{I\colon\;I\subset\{1,\dots,d\},\,I\neq\emptyset}
 W^x_I(\alpha,4t)\right\}\) \\
&\qquad \ge 2^{-d}e^{-d/2t^2}
\endaligned \tag4.3
$$
with the Gaussian probability measure $\mu_{n_1+\cdots+n_d,t}$ defined
above.}
\medskip
 
The following corollary of Proposition~4.2 is important for us.
 
\medskip\noindent
{\bf Corollary of Proposition~4.2.} {\it Let us have
a pseudometric $\rho_\alpha$ in $R^{n_1}\times\cdots\times R^{n_d}$
induced by a pseudonorm $\alpha$ in the tensor product
$R^{n_1}\otimes\cdots\otimes R^{n_d}$ of the Euclidean spaces
$R^{n_j}$, $1\le j\le d$. Let
$D\subset B^{n_1}\times\cdots\times B^{n_d}$ be a subset of the
product of the unit balls $B^{n_j}$, $1\le j\le d$ that has the
following property:
$\summ_{I\subset\{1,\dots,d\},\,I\neq\emptyset}W^x_I(\alpha,4t)\le u$
with some fixed numbers $0<t\le1$ and $u>0$ for all $x\in D$.
 
Then there is a constant $C>0$ depending only on the parameter $d$
such that the set $D$ has a $2u$-net of cardinality $e^{C/t^2}$
with respect to the pseudometric $\rho_\alpha$. In more detail this
means that there is a set $\{x^{(1)},\dots,x^{(N)}\}\subset D$ with
cardinality $N\le e^{C/t^2}$ such that
$\min\limits_{1\le j\le N}\rho_\alpha(x,x^{(j)})\le 2u$ for all
$x\in D$.}
 
\medskip\noindent
{\it Proof of the Corollary.} Let us construct a sequence
$x^{(1)},x^{(2)},\dots,x^{(N)}$, $x^{(j)}\in D$, $1\le j\le N$, in
the following way. Let us choose first a point $x^{(1)}\in D$ in an
arbitrary way. If the points $x^{(1)},\dots,x^{(j)}$ are already
chosen, and there are some points $x\in D$ such that
$\rho_\alpha(x,x^{(p)})>2u$ for all $1\le p\le j$, then we choose
an arbitrary point $x\in D$ with this property as $x^{(j+1)}$. If
there is no such point, then we finish our procedure at the $j$-th
step. Let $N$ be the number of points $x^{(j)}$ that we could
choose in such a way. Observe that the sets
$U_j=\{y\colon\;y\in R^{n_1}\times\cdots\times R^{n_d},
\;\rho_\alpha(y,x^{(j)})\le u\}$, $1\le j\le N$, are disjoint, 
because $\rho_\alpha(x_j,x_{j'})>2u$ for all $1\le j,j'\le N$,
$j\neq j'$. Beside this,
$\mu_{n_1+\cdots+n_d,t}(U_j)\ge 2^{-d}e^{-d/2t^2}$
by Proposition~4.2 for all $1\le j\le N$. Hence
$N\le 2^de^{d/2t^2}\le e^{C/t^2}$. Beside this, the set
$\{x^{(1)},\dots,x^{(N)}\}$ is a $2u$-net in $D$, because if 
there were a point $x\in D$ such that
$\min\limits_{1\le j\le N}\rho_\alpha(x,x^{(j)})>2u$ then we
would not finish our procedure at the $N$-th step.

\medskip\noindent
{\it Remark.} In the proof of the above corollary we applied a
rather standard method, well-known in the literature.
In general applications  of a result similar to
Proposition~4.2 the cardinality of a good $\e$-net of
the set $B^{n_1}\times\cdots\times B^{n_d}$ is bounded. Here
a slightly more general result was proved. This corollary gave
an estimate about the cardinality of a good $\e$-net of an arbitrary
set $D\subset B^{n_1}\times\cdots\times B^{n_d}$. For some sets
$D$ with nice properties it provides a much better bound for the
cardinality of a good $\e$-net in $D$ than for the cardinality of
a good $\e$-net in $B^{n_1}\times\cdots\times B^{n_d}$.
This observation will be exploited in our further considerations.
 
\medskip
In formula (2.11) we defined a pseudometric $\rho_\alpha$ in
the product $R^{n_1}\times\cdots\times R^{n_{d-1}}$ of the
Euclidean spaces $R^{n_j}$, $1\le j\le n$ and in formula~(2.12)
a pseudonorm $\alpha$ in their tensor product
$R^{n_1}\otimes\cdots\otimes R^{n_{d-1}}$. A comparison of
formulas~(2.13) and~(4.1) shows that Proposition~4.2 and its
corollary can be applied (with parameter $d-1$) for the metric
$\rho_\alpha$ and norm $\alpha$ defined in~(2.11) and~(2.12).
This fact plays an important role in the proof of the Basic
estimate.

\beginsection 5. The Main inequality.

In this section I formulate a result that I call the Main inequality
and show that the Basic estimate and in such a way Theorem~1 follows
from it. This result is a weaker version of an inductive statement
formulated in the proof of Theorem~3 in~[2]. I had to formulate such
a weaker statement because the corresponding result in~[2] seems to
be incorrect.
 
Let us fix the parameter $d\ge3$. We shall define appropriate classes
$\Cal U(r,N)$ depending on two parameters $N$ and $r$ which consist of
finite subsets of $R^{n_1}\times\cdots\times R^{n_{d-1}}$ with some
nice properties. In the Main inequality we give an estimate on the 
moments of the random variables
$\sup\limits_{u\in U,\,u'\in U}[Y_d(u)-Y_d(u')]$ for the sets 
$U\in\Cal U(r,N)$, where $Y_d(u)$  with parameter
$u\in R^{n_1}\times\cdots\times R^{n_{d-1}}$ is the Gaussian random 
variable defined in~(2.1). To define these classes of sets 
$\Cal U(r,N)$ some additional quantities have to be introduced.
 
We shall work with the linear functional $A(v)=A(v,d)$ defined for
functions $v\in R^{n_1}\otimes\cdots\otimes R^{n_d}$ in
formula~(1.3) with the help of a set of numbers
$A=\{a(i_1,\dots,i_d),\;1\le i_p\le n_p,\,1\le p\le d\}$. Let us
also recall the definition of the Gaussian random variables $Y_d(u)$
defined in~(2.1) for vectors
$u=(u_1,\dots,u_{d-1})\in R^{n_1}\times\cdots\times R^{n_{d-1}}$
together with a standard Gaussian random vector
$G_d=(g_d(1),\dots,g_d(n_d))$. We shall also work with the quantity
$\rho_\alpha(u,v)$, $u\in R^{n_1}\times\cdots\times R^{n_{d-1}}$
and $v\in R^{n_1}\times\cdots\times R^{n_{d-1}}$ defined in~(2.11).
 
Beside this, to define the sets $\Cal U(r,N)$ we still have to
introduce some pseudonorms $\tilde\alpha_{j,k}$ in the spaces 
$R^{n_j}$ for all pairs $j,k$ such that $1\le j,k\le d-1$, 
$j\neq k$, with the help of the coefficients $a(i_1,\dots,i_d)$ 
appearing in formula~(1.3).
 
For this goal first we introduce the set of constants
$$
\aligned
b_{u_j}^{(j)}(i_1,\dots,i_{j-1},i_{j+1},\dots,i_d)
&=\summ_{i_j\colon\;1\le i_j\le n_j}a(i_1,\dots,i_d)u_j(i_j),\\
&\qquad\qquad  1\le i_p\le n_p,  \; p\in\{1,\dots,d\}\setminus\{j\},
\endaligned \tag5.1
$$
for all vectors $u_j\in R^{n_j}$ and the functional
$$
B_{u_j}^{(j)}(v)=\sum \Sb (i_1,\dots,i_{j-1},i_{j+1},\dots,i_d) \\
1\le i_p\le n_p,\;p\in\{1,\dots,d\}\setminus\{j\} \endSb
b_{u_j}^{(j)}(i_1,\dots,i_{j-1},i_{j+1},\dots,i_d)
v(i_1,\dots,i_{j-1},i_{j+1},\dots,i_d) \tag5.2
$$
depending on this $u_j\in R^{n_j}$ for all
$v\in R^{n_1}\otimes\cdots\otimes R^{n_{j-1}}\otimes
R^{n_{j+1}}\otimes\cdots\otimes R^{n_d}$. The functional
$B_{u_j}^{(j)}(v)$ defined in~(5.2) is a special case of the
operator $B_K(v)$ introduced in~(1.4) if we choose
$K=\{1,\dots,d\}\setminus\{j\}$ and the coefficient $b_K(\cdot)$
are chosen as the numbers $b_{u_j}^{(j)}(\cdot)$ introduced
in~(5.1). With such a choice we can introduce the quantity
$V(P,B^{(j)}_{u_j})=V(P,b_{u_j}^{(j)}(\cdot))$ for all partitions
$P$ of the set $\{1,\dots,d\}\setminus\{j\}$ as the quantity
$V(P,B_K)=V(P,b_K(\cdot))$
defined in~(1.6) with this choice $K=\{1,\dots,d\}\setminus\{j\}$
and $B_K(v)=B_{u_j}^{(j)}(v)$. Let $P_{j,k}$ denote the partition
$P_{j,k}=\{\{k,d\},\{l\},1\le l\le d-1, l\neq j,k\}$ of the set
$\{1,\dots,d\}\setminus\{j\}$, and define
$$
\tilde\alpha_{j,k}(u_j)=V(P_{j,k},B_{u_j}^{(j)}), \quad 1\le j,k\le d-1,
\quad k\neq j,  \quad u_j\in R^{n_j}. \tag5.3
$$
It is easy to check that $\tilde\alpha_{j,k}(u_j)$ is a pseudonorm
in $R^{n_j}$.
 
The expression $\tilde\alpha_{j,k}(u_j)$ can also be written as
$$
\aligned
\tilde\alpha_{j,k}(u_j)=
&\sup\Sb v_p(\cdot),\,p\in\{1,\dots,d-1\}\setminus\{j,k\},\,
v_{k,d}(\cdot,\cdot)\colon \\
\summ_{i_p}v_p^2(i_p)\le 1,
\; p\in\{1,\dots,d-1\}\setminus\{j,k\},\;
\summ_{i_k,i_d}v^2(i_k,i_d)\le1,\endSb
\;\sum_{i_1,\dots,i_d}a(i_1,\dots,i_d) u_j(i_j)v(i_k,i_d)  \\
&\qquad\qquad\qquad\qquad\qquad
\prod_{p\in\{1,\dots,d-1\}\setminus\{j,k\}}v_p(i_p)
\endaligned \tag5.4
$$
for any $u_j=(u_j(1),\dots,u_j(n_j))\in R^{n_j}$.
 
Given an operator $A(v)$ of order $d$, $d\ge3$, defined in
(1.3) and a positive integer~$M$ the following classes of sets
$\Cal U(r,N)=\Cal U_{A,M,d}(r,N)$ consisting of at most $r$ elements
$u\in R^{n_1}\times\cdots\otimes R^{n_{d-1}}$ will be introduced.
$$
\aligned
\Cal U(r,N)&=\Cal U_{A,M,d}(r,N) \\
&= \biggl\{U=\{(u^{(t)}=
(u^{(t)}_1,\dots,u^{(t)}_{d-1})\in R^{n_1}\times\cdots\times
R^{n_{d-1}},\;1\le t\le r'\}\colon\\
&\qquad 1\le r'\le r, \quad
\tilde\alpha_{j,k}(u^{(t)}_j)\le 2^{-N}M^{-(d-2)/2},
\text{ for all }1\le t\le r' \\
&\qquad\qquad \text{ and } 1\le j,k\le d-1,\; j\neq k, \\
&\qquad \rho_\alpha(u^{(t)},u^{(t')})\le2^{-2N}M^{-(d-1)/2}
\quad \text{for all } 1\le t,t'\le r',\\
&\qquad u^{(t)}\in B^{n_1}\times \cdots\times B^{n_{d-1}},
\quad \text{for all } 1\le t\le r',\\
&\qquad u^{(t)}-u^{(t')}\in B^{n_1}\times \cdots\times B^{n_{d-1}}
\quad \text{for all }
1\le t,t'\le r'\biggr\}
\endaligned\tag5.5
$$
with the above defined $\tilde\alpha_{j,k}$ and the quantity
$\rho_\alpha(\cdot,\cdot)$ introduced in~(2.11).
 
In the Main estimate we shall prove a moment estimate for the
supremum of some random variables determined with the help of the
sets $\Cal U(r,N)$. It holds under the condition
$$
\alpha_s=\alpha_s(A)\le M^{-(s-1)/2},\quad \text{for all }
1\le s\le d-1,  \tag5.6
$$
where the quantities $\alpha_s$ were defined in~(1.7).
 
\medskip\noindent
{\it Remark.} In Theorem~1A we imposed a similar but stronger
condition in formula~(1.10). It also contained the condition
$\alpha_d\le M^{-(d-1)/2}$ for $s=d$. This condition is missing
here. It is replaced by the inequalities imposed on $\rho_\alpha$
in the definition of the sets $\Cal U(r,N)$. The additional
condition of Theorem~1A is needed when we want to prove the
Basic estimate with the help of the Main inequality.
 
\medskip\noindent
{\bf The Main inequality.} {\it Let a multilinear form $A$ of order
$d\ge3$ satisfy condition~(5.6). Take a standard normal random
vector  $G_d=(g_d(1),\dots,g_d(n_d))$ of dimension $n_d$, and
introduce with its help the random variables $Y_d(u)$ defined
in~(2.1) for all vectors $u=(u_1,\dots,u_{d-1})$,
$u_p=(u_p(1),\dots, u_p(n_p))\in R^{n_p}$, $1\le p\le d-1$.
There is a threshold index $A_0\ge0$ and a constant $C>0$ such
that for integers $r\ge1$ and $N\ge0$ the inequality
$$
E\[\sup\Sb u^{(t)}=(u_1^{(t)},\dots, u_{d-1}^{(t)})\in U,\\
u^{(t')}=(u_1^{(t')},\dots, u_{d-1}^{(t')})\in U \endSb
(Y_d(u^{(t)})-Y_d(u^{(t')})) \]^{2^{2(N+A)}M} \!\!\!\!
\le(CM^{-(d-2)/2}2^{(A-N)})^{2^{2(N+A)}M} \tag5.7
$$
holds for all $U\in\Cal U(r,N)$ and integers $A\ge A_0$. The
numbers $A_0=A_0(d)$ and $C=C(d)$ are sufficiently large
constants which depend only on $d$ and do not depend on the
parameters $r$ and~$N$.}
 
\medskip
Now I give the proof of the Basic estimate with the help of the
Main inequality.
 
\medskip\noindent
{\it The proof of the Basic estimate.}\/ First we show that under
the conditions of the Basic estimate $U\in \Cal U(r,0)$ for any set
$U=\{(u^{(t)},\;1\le t\le r\}$ consisting of $r$ vectors
$u^{(t)}=(u_1^{(t)},\dots,u^{(t)}_{d-1})$, $1\le t\le r$ such that
$2u_j^{(t)}\in B^{n_j}$ for all $1\le j\le d-1$, $1\le t\le r$.
To show this observe that
$$
\rho_\alpha(u^{(t)},u^{(t')})\le\rho_\alpha(u^{(t)},0)+
\rho_\alpha(u^{(t')},0)
$$
for all $1\le t,t'\le r$, where $0$ denotes the vector with all
coordinates $0$, and
$$
\align
\rho_\alpha(u^{(t)},0)
&\le 2^{-(d-1)}\sup\Sb u=(u_1,\dots,u_{d-1})\colon\\
u_j\in B^{n_j},\;1\le j\le d-1 \endSb \rho_\alpha(u,0) \\
&= 2^{-(d-1)}\sup\Sb u=(u_1,\dots,u_{d-1})\colon\\
u_j\in B^{n_j},\;1\le j\le d-1 \endSb
\[\sum_{i_d}\(\sum_{i_1,\dots,i_{d-1}}
a(i_1,\dots,i_d)\prod_{j=1}^{d-1} u_j(i_j)\)^2\]^{1/2}\\
&= 2^{-(d-1)}\sup\Sb u=(u_1,\dots,u_d)\colon\\
u_j\in B^{n_j},\;1\le j\le d \endSb
\sum_{i_1,\dots,i_d} a(i_1,\dots,i_d)
\prod_{j=1}^d u_j(i_j)= 2^{-(d-1)}\alpha_d \\
&\le 2^{-(d-1)}M^{-(d-1)/2}
\endalign
$$
for all $1\le t\le r$, and a similar estimate holds for
$\rho_\alpha(u^{(t')},0)$. (This is the point where we exploited
that the estimate $\alpha_s\le M^{-(s-1)/2}$ also holds for $s=d$.)
Beside this
$\alpha_{j,k}(u_j^{(t)})\le\frac12\alpha_{d-1}\le\frac12M^{-(d-2)/2}$,
and clearly $u^{(t)}\in B^{n_1}\times\cdots\times B^{n_d}$,
and $u^{(t)}-u^{(t')}\in B^{n_1}\times\cdots\times B^{n_d}$ for
all $1\le t,t'\le r$. The above relations imply that
$U\in\Cal U(r,0)$.
 
It can be proved with the help of the above fact and the Main
inequality with the choice $N=0$ that
$$
E\[\sup_{u\in\frac12B^{n_1}\times\cdots\times
\frac12B^{n_{d-1}}}Y_d(u)\]^{2^{2A_0}M}
\le \(C2^{2A_0}M^{(d-2)/2}\)^{2^{A_0}M} \tag5.8
$$
with the same number $A_0$ which appears in the Main inequality as
the threshold index.
 
To prove this statement let us list the set of vectors
$u\in \frac12 B^{n_1}\times\cdots\times\frac12B^{n_{d-1}}$ such
that all their coordinates are rational numbers in a sequence $u^{(t)}$,
$t=1,2,\dots$. Let $u^{(1)}=(0,\dots,0)$ in this sequence. Let
$U_r=\{u^{(t)},\;1\le t\le r\}$ be the set consisting of the
first $r$ terms of this sequence. Observe that
$$
\sup_{u\in B^{n_1}\times\cdots\times B^{n_{d-1}}} Y_d(u)
=\lim_{r\to\infty}\sup_ {u^{(t)}\in U_r}Y_d(u^{(t)})
$$
Let us apply a weakened form of the Main inequality with $N=0$ and
$A=A_0$ (we may assume that $A_0\ge1$) for all above defined
sets $U_r$, $r=1,2,\dots$, where instead of taking the supremum of all
differences $Y_d(u^{(t)})-Y_d(u^{(t')})$, $1\le t,t'\le r$ we take
this supremum only for pairs $(t,t')$ with $t'=1$. In this case
$Y_d(u^{(t')})=0$ with probability 1. The series of inequalities
obtained in such a way, (where the upper bound does not depend
on~$r$) together with the previous identity and the Beppo-Levy
theorem imply relation~(5.8). This inequality together with the
H\"older inequality for $p=2^{2A_0-1}$ yield that
$$
E\[\sup_{u\in\frac12B^{n_1}\times\cdots\times
\frac12B^{n_{d-1}}}Y_d(u)\]^{2M}
\le \(C2^{A_0}M^{(d-2)/2}\)^{2M} \le\bar C^MM^{-(d-2)M} \tag5.9
$$
with a universal constant $\bar C$. Relation (2.6) follows from
this inequality. To see this it is enough to observe that
if the condition
$u\in\frac12B^{n_1}\times\cdots\times\frac12B^{n_{d-1}}$ is replaced
by the condition $u\in B^{n_1}\times\cdots\times B^{n_{d-1}}$,
then the inequality remains valid if the right-hand side is
multiplied by $2^{(d-1)M}$, i.e. the constant $\bar C$ is multiplied
by $2^{(d-1)}$ in~(5.9).
 
\medskip\noindent
{\it Remark.} Actually we needed the Main inequality only for $N=0$
(and arbitrary $r$). But we shall prove it by a backward induction
procedure. It is not difficult to see that the Main inequality
holds if $N\ge N_0$ with a very large $N_0$ whose value may depend
on~$r$. If this is shown, then we may apply
backward induction to prove the Main inequality. It may seem a
technical point that the hardest estimate of this paper is proved
by a backward and not by a forward induction. But I think that the
situation is much more complex.
 
I met a similar situation in a study leading to paper~[5]. Here
also backward induction had to be applied to solve the hardest
part of the problem, and this had a non-technical reason. The
supremum of such random variables had to be bounded whose
behaviour was very `non-Gaussian'. The main contribution to the
supremum I was interested in came from the influence of some
irregular events. These irregular events had very small
probability, but they played a dominant role because of their
large number. Their effect could be controlled by means of a
backward and not by a forward induction procedure. I believe that 
behind the proof of the Main inequality in this paper a similar 
phenomenon is hiding. But to understand the situation better some 
additional work has to be done.
 
\beginsection 6. Some results about the existence of good partitions.
 
The proof of the Main inequality is based on the existence of
some good partitions of the class of sets $\Cal U(r,N)$ defined
in~(5.5). These results are formulated in this section in Lemmas~6.1
and~6.2. Lemma~6.1 can be considered as a version of Lemma~8 in~[2],
and Lemma~6.2 is an improvement of this result. It states that there
exists a partition of the sets $U\in\Cal U(r,N)$ which satisfies
Lemma~6.1, and it also has some extra properties useful in our
investigation. Its cardinality can be bounded similarly to
Lemma~6.1. Such a result was needed to get a proof without the
application of Theorem~3 of~[2] whose validity is questionable.
Lemmas~6.1 and~6.2 are proved by means of Proposition~4.2 and
its corollary. But to prove them we also need some additional
inequalities. They are given in Lemma~6.3 which can be
considered as a version of Lemmas~5 and~6 in~[2]. Lemma~6.3 is
formulated and proved in this section.
 
Before the formulation of these results some additional notations
have to be introduced. We define with the help of a vector
$u\in R^{n_1}\times\cdots\times R^{n_{d-1}}$ and a set
$I\subset\{1,\dots,d-1\}$ an operator which is a special case of
the class of operators defined in formula~(1.4). We also
introduce some quantities corresponding to this operator which
are the analogs of the quantities $\alpha_s$,
$\tilde\alpha_{j,k}$, $\rho_\alpha$ defined earlier with the
help of the operator~$A(v)$ given in~(1.3).
 
Fix a set $I=\{j_1,\dots,j_s\}\subset\{1,\dots,d-1\}$
with $1\le s\le d-2$ elements and a vector
$u=(u_1,\dots,u_{d-1})\in R^{n_1}\times\cdots\times R^{n_{d-1}}$.
Let us define with their help the numbers
$$
b_u^{I}(i_j,\;j\in \{1,\dots,d\}\setminus I)
=M^{|I|/2}\summ_{(i_j,\; j\in I)}
a(i_1,\dots,i_d)\prod_{j\in I}u_j(i_j) \tag6.1
$$
depending on the vectors $(i_j,\; j\in\{1\dots,d\}\setminus I)$
and the linear functional
$$
B_u^{I}(v)=\sum_{(i_j,\;j\in \{1,\dots,d)\setminus I)}
b_u^{I}(i_j,\;j\in \{1,\dots,d\}\setminus I)
v(i_{j},\;  j\in \{1,\dots,d\}\setminus I) \tag6.2
$$
acting on the space of functions $v=v(i_{j_1},\dots,i_{j_p})
\in R^{n_{j_1}}\otimes\cdots\otimes R^{n_{j_p}}$, with the set of
indices $\{j_1,\dots,j_p\}=\{1,\dots,d\}\setminus I$.
 
This operator
$B_u^I(v)$ is a special case of the operators $B_K(v)$ defined in
formula~(1.4) when $K=\{1,\dots,d\}\setminus I$, and coefficients
$b_K(\cdot)$ are the numbers $b_u^I(\cdot)$ defined in~(6.1).
(In the definition of the coefficients
$b_u^{I}(i_j,\;j\in \{1,\dots,d\}\setminus I)$ in~(6.1) a
multiplying factor $M^{|I|/2}$ was inserted. I applied such a
norming factor, because it simplifies the subsequent calculations.)
 
We can define the quantities $V(P,B_u^I)$ for all partitions $P$ of
the set $K=\{1,\dots,d\}\setminus I$ by formula~(1.6) with the choice
$B_K(v)=B_u^I(v)$. Let us also introduce, similarly to $\alpha_s$
defined in~(1.7) the quantity
$$
\alpha_{u,s}(I)=\supp_{P\colon\; |P|=s}V(P,B_u^{I}), \tag6.3
$$
where all partitions $P$ of the set $\{1,\dots,d\}\setminus I$ with 
cardinality $s$ are taking part in the supremum. We also introduce 
the numbers
$$
\tilde\alpha^I_{k}(u)=V(P_{I,k},B^I_u) \quad \text{for all }
k\in\{1,\dots,d-1\}\setminus I \tag6.4
$$
with the help of formula~(1.6), where the operator $B^I_u(v)$
defined in~(6.2) plays the role of $B_K(v)$, and the partition
$P_{I,k}$ of the set $\{1,\dots,d\}\setminus I$ is defined as
$P_{I,k}=\{\{k,d\},\{l\},\;l\in\{1,\dots,d-1\}\setminus(I\cup\{k\})\}$.
 
We shall also work with a quantity $\rho_{\alpha^I_u}(v,\bar v)$
defined for all pairs $(v,\bar v)$,
$v\in R^{n_1}\times\cdots\times R^{n_{d-1}}$ and
$\bar v\in R^{n_1}\times\cdots\times R^{n_{d-1}}$ with the help of
a vector $u\in R^{n_1}\times\dots\times R^{n_{d-1}}$ and set
$I\subset\{1,\dots,d-1\}$, $1\le |I|\le d-2$, similarly to the
term $\rho_\alpha$ introduced in~(2.11). First we define a version
of it.
$$
\aligned
\bar\rho_{\alpha^I_u}(v,\bar v)
=\biggl(&\sum_{i_d}\biggl[\sum_{(i_j,\;j\in\{1,\dots,d-1\}\setminus I)}
b^I_u(i_j,\;j\in\{1,\dots,d\}\setminus I)\\
&\qquad\biggl( \prod_{j\in\{1,\dots,d-1\}\setminus I} v(i_j)
-\prod_{j\in\{1,\dots,d-1\}\setminus I}\bar v(i_j)\biggr)\biggr]^2\biggl)^{1/2}
\endaligned \tag6.5
$$
for pairs of vectors
$v=(v_{j_1},\dots,v_{j_p})\in R^{n_{j_1}}\times\cdots\times R^{n_{j_p}}$
and $\bar v=(\bar v_{j_1},\dots,\bar v_{j_p})
\in R^{n_{j_1}}\times \cdots\times R^{n_{j_p}}$, where
$\{j_1,\dots,j_p\}=\{1,\dots,d-1\}\setminus I$. Observe that
$\bar\rho_{\alpha^I_u}$ is the pseudometric induced by the pseudonorm
$$
\aligned
\alpha^I_u(v)
&=\alpha^I_u(v(i_j, j\in\{1,\dots,d-1\}\setminus I))\\
&=\biggl(\sum_{i_d}\biggl[\sum_{(i_j,\;j\in\{1,\dots,d-1\}\setminus I)}
b^I_u(i_j,\;j\in\{1,\dots,d\}\setminus I)\\
&\qquad\qquad\qquad
v(i_j,\;j\in\{1,\dots,d-1\}\setminus I)\biggr]^2\biggl)^{1/2}
\endaligned \tag6.6
$$
on the tensor product $R^{n_{j_1}}\otimes \cdots\otimes R^{n_{j_p}}$,
where $\{j_1,\dots,j_p\}=\{1,\dots,d-1\}\setminus I$.
 
We can define the metric $\rho_{\alpha^I_u}$ in the space
$R^{n_1}\times\cdots\times R^{n_{d-1}}$ with the help of the metric
$\bar\rho_{\alpha^I_u}$ defined in~(6.5). To do this we introduce
the following notation. Given a vector
$v=(v_1,\dots,v_{d-1})\in R^{n_1}\times \cdots\times R^{n_{d-1}}$
and a set $I\subset\{1,\dots,d-1\}$ let $v_{I^c}$ denote the vector
we obtain by omitting the coordinates of the vector~$v$ belonging
to the set~$I$, i.e. let
$v_{I^c}\in R^{n_{j_1}}\times\cdots\times R^{n_{j_p}}$, and
$v_{I^c}=(v_j,\; j\in\{1,\dots,d-1\}\setminus I)$. Given two vectors
$v=(v_1,\dots,v_{d-1})\in R^{n_1}\times \cdots\times R^{n_{d-1}}$
and $\bar v=(\bar v_1,\dots,\bar v_{d-1})
\in R^{n_1}\times \cdots\times R^{n_{d-1}}$ put
$$
\rho_{\alpha^I_u}(v,\bar v)=\bar\rho_{\alpha^I_u}(v_{I^c},\bar v_{I^c}).
\tag6.7
$$
 
Now I formulate Lemma~6.1 and its strengthened version Lemma~6.2.
 
\medskip\noindent
{\bf Lemma 6.1.} {\it If an operator $A$ of order $d\ge3$ satisfies
relation (5.6), then each set $U\in \Cal U(r,N)=\Cal U_{A,M,d}(r,N)$
has a partition $u^{(1)}+U_1$, $u^{(2)}+U_2$,\dots, $u^{(L)}+U_L$
with $L\le 2^{C(d)M 2^{2N}}$ elements such that $U_l\in\Cal U(r,N+2)$
and $u^{(l)}\in U$ for all $1\le l\le L$. The number $C(d)$ depends
only on the order $d$ of the operator~$A$.}
 
\medskip\noindent
{\bf Lemma 6.2.} {\it Under the conditions of Lemma 6.1 each set
$U\in \Cal U(r,N)$ has a partition $u^{(1)}+U_1,u^{(2)}+U_2,\dots,
u^{(L)}+U_L$ with $L\le 2^{C(d)M 2^{2N}}$ elements such that
$u^{(l)}\in U$, $U_l\in\Cal U(r,N+2)$, $1\le l\le L$, and it
also satisfies following additional property. The inequality
$\rho_{\alpha_{u^{(l)}}^I}(u,\bar u)\le2^{-2N}M^{-(d-|I|-1)/2}$
holds for all sets $I\subset\{1,\dots,d-1\}$, $1\le |I|<d-2$, and
pairs of elements $u\in U_l$ and $\bar u\in U_l$, $1\le l\le L$.
The vector $u^{(l)}$ in this inequality is the same vector which
appears in the definition of the element $u^{(l)}+U_l$ of the
partition of $U$. The quantity $\rho_{\alpha^I_u}(\cdot,\cdot)$ was
defined in~(6.5) and~(6.7).}
 
\medskip
Lemmas 6.1 and 6.2 will be proved with the help of the following
Lemma~6.3.
 
\medskip\noindent
{\bf Lemma 6.3.} {\it Let a functional $A(v)$ of order $d\ge3$
defined in (1.3) satisfy condition~(5.6). Then for any
$u=(u_1,\dots,u_{d-1})\in B^{n_1}\times\cdots\times B^{n_{d-1}}$
the quantities $W^u_I(\alpha,t)$, $I\in\{1,\dots,d-1\}$,
$I\neq\emptyset$, defined in~(4.2) with the pseudonorm $\alpha$
introduced in~(2.12) satisfy the following inequalities.
$$
W^u_I(\alpha,t)\le\frac{t^{|I|}}{M^{(d-|I|-1)/2}} \quad
\text{if}\quad 2\le|I|\le d-1. \tag6.8
$$
For a set $I=\{k\}$ containing one element
$$
W^u_{\{k\}}(\alpha,t)\le t\min_{1\le j\le d-1,\, j\neq k}
\tilde\alpha_{j,k}(u_j), \tag6.9
$$
where $\tilde\alpha_{j,k}(u_j)$ was defined in~(5.3). Beside this,
$$
E\tilde\alpha_{j,k}(G_j)\le\frac {C(d)}{M^{(d-3)/2}}
\quad \text{for all }1\le j,k\le d-1,\; j\neq k, \tag6.10
$$
where $C(d)$ depends only on $d$, and $G_j$ is a 
standard normal vector of dimension~$n_j$.}
 
\medskip\noindent
{\it The proof of Lemma~6.3.}\/ For any set $I\subset\{1,\dots,d-1\}$,
$I\neq\emptyset$ and $u\in B^{n_1}\times\cdots\times B^{n_{d-1}}$
$$
\align
W^u_I(\alpha,1)&=E\(\[\sum_{i_d}\(\sum_{i_1,\dots,i_{d-1}} a(i_1,\dots,i_d)
\prod_{j\in\{1,\dots,d-1\}\setminus I} u_j(i_j)
\prod_{j\in I} g_j(i_j)\)^2\]\)^{1/2} \\
&\le \(E\[\sum_{i_d}\(\sum_{i_1,\dots,i_{d-1}} a(i_1,\dots,i_d)
\prod_{j\in\{1,\dots,d-1\}\setminus I} u_j(i_j)
\prod_{j\in I} g_j(i_j)\)^2\]\)^{1/2}\\
&=\[\sum_{(i_p,\;p\in I\cup\{d\})}
\(\sum_{(i_j,\; j\in \{1,\dots,d-1\}\setminus I)} a(i_1,\dots,i_d)
\prod_{j\in\{1,\dots,d-1\}\setminus I} u_j(i_j)\)^2\]^{1/2}\\
&=\sup \Sb v(i_p,\;p\in I\cup\{d\}) \colon\\
\sum v^2(i_p,\;p\in I\cup\{d\})\le1 \endSb
\sum_{i_1,\dots,i_d} a(i_1,\dots,i_d)
v(i_p,\;p\in I\cup \{d\})
\prod _{j\in\{1,\dots,d-1\}\setminus I}u_j(i_j) \\
&\le V(P_I,A), \tag6.11
\endalign
$$
where $P_I$ is the partition $P_I=\{I\cup\{d\},
\{j\}, 1\le j\le d-1,\, j\notin I\}$ of the set $\{1,\dots,d\}$, and
$V(P,A)$ is defined in~(1.6).
 
Since the partition $P_I$ has $d-|I|$ elements this inequality
together with relation (5.6) imply that for $|I|\ge2$
$$
W^u_I(\alpha,t)=t^{|I|}W^u_I(\alpha,1)\le t^{|I|}\alpha_{P_I}(A)
\le\frac{t^{|I|}}{M^{(d-|I|-1)/2}},
$$
i.e. (6.8) holds. In the case $I=\{k\}$ we get from the last
but one bound in (6.11), the representation of
$\tilde\alpha_{j,k}(u_j)$ in formula~(5.4) and the choice of an
arbitrary point $j\in\{1,\dots,d-1\}\setminus\{k\}$ that
$W^u_I(\alpha,1)\le\tilde\alpha_{j,k}(u_j)$, and this relation
implies formula~(6.9).
 
Inequality~(6.10) can be deduced from inequality (2.6) in the Basic
estimate with parameter $d-1$ if we write up the expression
$\tilde\alpha_{j,k}(G_j)$ in the form~(5.4), (by replacing the
vector $u_j$ by $G_j$ in it), consider it as an expression of the
form~(2.1) with $d-1$ variables by taking the pair $(k,d)$ as one
variable. Let us observe that relation~(5.6) implies 
relation~(1.10) with parameter $d-1$ in this case, hence we may
apply the Basic estimate. Let us apply a reindexation of the arguments 
by which the $j$-th variable turns to the $d-1$-th coordinate. The 
Basic estimate remains valid after such a reindexation. Since in the 
proof of Lemma~6.3 for parameter~$d$ we may assume that the Basic 
estimate holds for~$d-1$ we get inequality~(6.10) from the Basic 
estimate and the estimate $EZ_d\le(EZ_d^{2M})^{1/2M}$ which is a 
consequence of H\"older's inequality.
 
\beginsection 7. The proof of Lemmas 6.1 and 6.2 about the
existence of good partitions.
 
In this section Lemmas 6.1 and 6.2 will be proved with the help
of Proposition~4.2, its corollary and Lemma~6.3.
 
\medskip\noindent
{\it The proof of Lemma~6.1.}\/ If relation (5.6) holds, then
relation~(6.10) in Lemma~6.3 implies the inequality
$E\tilde\alpha_{j,k}(G_j)\le CM^{-(d-3)/2}$ for all
$1\le j,k\le d-1$, $j\neq k$. Hence Proposition~4.1 yields
the estimate 
$\mu_{n_j,t}(y\colon\;y\in R^{n_j},\tilde\alpha_{j,k}(x-y)\le CtM^{-(d-3}/2)
\ge e^{-C'/t^2M^{(d-3)}}$ 
for all numbers $t>0$, pairs $(j,k)$,
$1\le j,k\le d-1$, $j\neq k$, and $x\in B^{n_j}$, where
$\mu_{n_j,t}$ denotes the distribution of $tG_j$ if $G_j$ is a
standard normal vector of dimension $n_j$. This estimate, or in 
a simpler way corollary of Proposition~4.2 yields the following 
result for the metric $\rho_\alpha(x,y)=\tilde\alpha_{j,k}(x-y)$ 
in the space $R^{n_j}$ with the choice $D=B^{n_j}$, 
$t=C2^{-N}M^{-1/2}$ and $u=2^{-(N+3)}M^{-(d-2)/2}$ 

For all pairs $(j,k)$, $1\le j,k\le d-1$, $j\neq k$, the unit ball 
$B^{n_j}\subset R^{n_j}$ has a partition
$\hat U^{(j,k)}_1,\dots,\hat U^{(j,k)}_{L(j,k)}$ with
$L(j,k)\le e^{C/t^2}\le2^{C2^{2N}M}$ elements such that
$\tilde\alpha_{j,k}(y-x)\le2^{-(N+2)}M^{-(d-2)/2}$ if 
$x\in\hat U_l^{(j,k)}$ and $y\in\hat U_l^{(j,k)}$ with the same 
index~$l$. Hence any set 
$U\subset B^{n_1}\times\cdots\times B^{n_{d-1}}$, in particular any 
set $U\in\Cal U(r,N)$ has a partition
$U_1^{(j,k)},\dots,U_{L(j,k)}^{(j,k)}$ with
$L(j,k)\le 2^{C2^{2N}M}$ elements  such that
$\tilde\alpha_{j,k}(x_j-y_j)\le 2^{-(N+2)}M^{-(d-2)/2}$ if
$x=(x_1,\dots,x_{d-1})\in U_l^{(j,k)}$ and
$y=(y_1,\dots,y_{d-1})\in U_l^{(j,k)}$ with the same index
$1\le l\le L(j,k)$.  Indeed, the sets
$U_l^{(j,k)}=\{y=(y_1,\dots,y_{d-1})\colon\;
y\in U,\,y_j\in\hat U_l^{(j,k)}\}$, $1\le l\le L(j,k)$, provide
such a partition of~$U$.
 
I claim that the existence of such partitions for all pairs
$(j,k)$, $1\le j,k\le d-1$, $j\neq k$, implies that each set
$U\in \Cal U(r,N)$ has a partition of the form
$u^{(1)}+\bar U_1$, $u^{(2)}+\bar U_2$,\dots,
$u^{L}+\bar U_{L}$ with $L\le 2^{C(d)M 2^{2N}}$ elements such
that $u^{(l)}\in U$, and
$\tilde\alpha_{j,k}(u_j)\le 2^{-(N+2)}M^{-(d-2)/2}$ if
$u=(u_1,\dots,u_{d-1})\in\bar U_l$ with some $1\le l\le L$
for all $1\le j,k\le d-1$, $j\neq k$.
 
To show this let us consider for all pairs $(j,k)$,
$1\le j,k\le d-1$, $j\neq k$, a partition
$U_1^{(j,k)},\dots,U_{L(j,k)}^{(j,k)}$ of the set $U$ with
$L(j,k)\le 2^{C2^{2N}M}$ elements such that
$\tilde\alpha_{j,k}(x_j-y_j)\le2t$ if
$x=(x_1,\dots,x_{d-1})\in U_l^{(j,k)}$ and
$y=(y_1,\dots,y_{d-1})\in U_l^{(j,k)}$ with the same index
$l=l(j,k)$. Take all intersections of the form
$\bigcap\limits_{(j,k)\colon\;1\le j,k\le d-1,j\neq k}U_{l(j,k)}^{(j,k)}$,
i.e. take all possible intersections which contain exactly
one element from each of the above partitions indexed by the pairs
$(j,k)$. By reindexing the sets obtained in such a way
we get a partition $\tilde U_1,\dots,\tilde U_L$ of the set $U$
with $L\le 2^{C2^{2N}M}$ elements such that for all pairs
$u=(u_1,\dots,u_{d-1})\in \tilde U_l$ and
$\bar u=(\bar u_1,\dots,\bar u_{d-1})\in \tilde U_l$ with the same
index $l$ and $1\le j,k\le d-1$, $j\neq k$,
$\tilde\alpha_{j,k}(u_j-\bar u_j)\le2^{-(N+2)}M^{-(d-2)/2}$. Then
choosing an arbitrary element $u^{(l)}\in \tilde U_l$ and writing
$\tilde U_l=u^{(l)}+\bar U_l$ with
$\bar U_l=\{u-u^{(l)}\colon\; u\in\tilde U_l\}$ we get a partition
with the desired property.
 
It can be shown with the help of the corollary of Proposition~4.2 
with the choice that each set $\bar U_l$, taking part in the 
above constructed partition $u^{(l)}+\bar U_l$, $1\le l\le L$, 
of the set $U$ has a partition $U_{l,1}\dots, U_{l,L_l}$ with 
$L_l\le 2^{C2^{2N}M}$ elements such that
$\rho_\alpha(u,\bar u)\le2^{-2(N+2)}M^{-(d-1)/2}$ if $u\in U_{l,p}$
and $\bar u\in U_{l,p}$ with the same parameters $l$ and $p$.
Indeed, let us choose $t=c2^{-N}M^{-1/2}$ with a suffficiently 
small constant $1\ge c>0$. Observe that with the choice of such a 
number $t$ and a vector $u\in\bar U_l$ with some index 
$1\le l\le L$ we can write by (6.8)
$$
W^u_I(\alpha,t)\le\frac{t^{|I|}}{M^{(d-|I|-1)/2}}\le c^2 2^{-2N}M^{-(d-1)/2}
$$
for all sets $I\subset\{1,\dots,d-1\}$ such that $|I|\ge2$.
For a set $I=\{k\}$, $1\le k\le d-1$, containing one element we have
$$
W^u_{\{k\}}(\alpha,t)\le t\min_{1\le j\le d-1,\, j\neq k}
\tilde\alpha_{j,k}(u_j)\le c2^{-2N}M^{-(d-1)/2}
$$
by relations~(6.9) and
$\tilde\alpha_{j,k}(u_j)\le 2^{-(N+2)}M^{-(d-2)/2}$ if
$u=(u_1,\dots,u_{d-1})\in\bar U_l$. Hence
$$
\sum_{I\colon\;I\subset\{1,\dots,d-1\},\,I\neq\emptyset}
W^u_I(\alpha,4t)\le 2^{-2(N+3)}M^{-(d-1)/2}
$$
for a vector $u\in\bar U_l$ if the parameter $c>0$ is chosen
sufficiently small. Then an application of the corollary of
Proposition~4.2 for one of the sets $\bar U_l$, $1\le l\le L$ with
the metric $\rho_\alpha$ and the choice $t=c2^{-N}M^{-1/2}$ and
$u=2^{-2(N+3)}M^{-(d-1)/2}$ shows that there exists a partition
$U_{l,p}$, $1\le p\le L_l$, of $\bar U_l$ of cardinality
$L_l\le 2^{C_1/t^2}\le 2^{CM2^{2N}}$ with the desired property.
 
Put $u^{(l,p)}=u^{(l)}$ for all $1\le l\le L$ and $1\le p\le L_l$,
and consider all sets $u^{(l,p)}+U_{l,p}$, $1\le l\le L$,
$1\le p\le L_l$. I claim that a reindexation of these sets provides
a partition of the set $U\in\Cal U(r,N)$ that satisfies Lemma~6.1.
Indeed, these sets provide a partition of the set $U$ with
$L\le 2^{CM2^{2N}}$ elements. Beside this, $u^{(l,p)}\in U$ for
all indices $l$ and $p$. We still have to check that
$U_{l,p}\in\Cal U(r,N+2)$ for all pairs of indices $l$ and $p$.
The elements of the sets $U_{l,p}$ satisfy the desired inequalities
for $\tilde\alpha_{j,k}$ and $\rho_\alpha$, and the sets $U_{l,p}$
have at most $r$ elements. To check that the sets $U_{l,p}$ satisfy
the remaining properties of the elements of the class
$\Cal U(r,N+2)$ observe that for a point $u\in U_{l,p}$
$u=\tilde u-u^{(l)}$ with $\tilde u\in \tilde U_{l,p}\subset U$ and
$u^{(l)}\in U$, hence $u\in B^{n_1}\times\cdots\times B^{n_{d-1}}$.
The analogous statement also holds for a difference $u-u'$ with
$u\in U_{l,p}$ and $u'\in U_{l,p}$, since such a difference can be
written as the difference of two vectors from the set
$\tilde U_l\subset U$.
 
\medskip\noindent
{\it The proof of Lemma 6.2.}\/ The main step of the proof is
the verification of the following statement formulated in
relation~(7.1).
 
Take a partition
$u^{(l)}+U_l$, $1\le l\le L$, of a set $U\in\Cal U(r,N)$ that 
satisfies Lemma~6.1, and fix one of the vectors $u^{(l)}$ in this 
partition together with a set $I\subset\{1,\dots,d-1\}$, 
$1\le |I|\le d-2$.  There is a partition 
$V_1=V_1(l,I),\dots,V_L=V_{L(l,I)}(l,I)$ with
$L(l,I)\le 2^{C2^{2N}M}$ elements of the product of unit
balls $B^{n_{j_1}}\times\cdots\times B^{n_{j_r}}$ with indices
$\{j_1,\dots,j_r\}=\{1,\dots,d-1\}\setminus I$ such that
$$
\bar\rho_{\alpha_{u^{(l)}}^I} (v,\bar v)\le2^{-2N}M^{-(d-|I|-1)/2}
\quad\text{if } v\in V_{p}(l,I)\text{ and }\bar v\in V_{p}(l,I)
\text { with an index } p,  
\tag7.1
$$
i.e. this inequality holds if $v$ and $\bar v$ are contained in the 
same element of the partition $V_p(l,I)$, $1\le p\le L(l,I)$, of the 
set $B^{n_{j_1}}\times\cdots\times B^{n_{j_r}}$. The metric
$\bar\rho_{\alpha_{u}^I}(v,\bar v)$ (with a general vector
$u\in R^{n_1}\times\cdots\times R^{n_{d-1}}$) was defined in 
formula~(6.5).

First the following inequalities will be verified. For all sets $I$,
$I\subset\{1,\dots,d-1\}$, $1\le |I|\le d-2$
$$
\alpha_{u^{(l)},s}(I)\le M^{-(s-1)/2} \quad\text{for all }
1\le s\le d-|I|-1 \tag7.2
$$
and
$$
\tilde\alpha^I_{k}(u^{(l)})\le 2^{-N}M^{-(d-|I|-2)/2} \quad
\text{for all } k\in\{1,\dots,d-1\}\setminus I, \tag7.3
$$
where $\alpha_{u,s}(I)$ was defined in~(6.3) and
$\tilde\alpha^I_k(u)$ in (6.4) (for a general vector $u$).
 
To check (7.2) let us compare a partition $P$ of
$\{1,\dots,d\}\setminus I$ of cardinality $|P|=s$, $1\le s\le d-|I|-1$,
with the partition $\bar P$ of the set $\{1,\dots,d\}$ we get by
attaching all one point sets of $I$ to the elements of the  partition $P$.
Then $|\bar P|=s+|I|$, hence $V(\bar P,A)\le \alpha_{s+|I|}(A)
\le M^{-(s+|I|-1)/2}$ by relation~(5.6) and
$V(P,B^I_{u^{(l)}})\le M^{|I|/2}V(\bar P,A)\le M^{-(s-1)/2}$.
Since this relation holds for all partitions $P$ such that
$|P|=s$ this implies~(7.2).
 
Beside this the relation
$u^{(l)}\in U$ with an $U\in\Cal U(r,N)$ implies that
$\tilde\alpha_{j,k}(u^{(l)}_j)\le 2^{-N}M^{-(d-2)/2}$, and
$\tilde\alpha^I_{k}(u^{(l)})\le M^{|I|/2} \tilde\alpha_{j,k}(u^{(l)}_j)
\le 2^{-N}M^{-(d-|I|-2)/2}$ for all $j\in I$ and
$k\in\{1,\dots,d-1\}\setminus I$. Hence relation~(7.3) also holds.
 
First we prove the existence of a partition with less than
$2^{C2^{2N}M}$ elements satisfying~(7.1) only in the case
$|I|\le d-3$. This will be done with the help of the corollary
of Proposition~4.2 when it is applied to the metric
$\bar\rho_{\alpha_{u^{(l)}}^I}$ and the norm $\alpha_{u^{(l)}}^I$
inducing it. These quantities were introduced in~(6.5) and~(6.6). 
In the proof we need good estimates on the terms
$W_K^u(\alpha_{u^{(l)}}^I,t)$ defined in~(4.2) for all sets 
$K\subset\{1,\dots,d-1\}\setminus I$, $K\neq\emptyset$ and
$u\in B^{n_{j_1}}\times\cdots\times B^{n_{j_s}}$ with a number~$t$
chosen as $t=c2^{-N}M^{-1/2}$ with a sufficiently small constant 
$1\ge c>0$. This quantity will be bounded by means of the 
estimates~(7.2), (7.3) and Lemma~6.3. More precisely an
equivalent version of Lemma 6.3  will be applied where
$B^I_{u^{(l)}}$ (defined in~(6.2)) is chosen as the operator $A$,
and as a consequence $\alpha_{u^{(l)},s}(I)$ defined in~(6.3) plays
the role of the term $\alpha_s=\alpha_s(A)$. This term must satisfy
relation~(5.6) to have the right to apply Lemma~6.3. (Actually the 
variables of the operator $B^I_{u^{(l)}}$ have to be reindexed if 
we want to apply Lemma~6.3 in its original form.) The operator 
$B^I_{u^{(l)}}$ acts on the functions on 
$\{1,\dots,d\}\setminus I$, on a set of $d-|I|$ elements, and by 
relation (7.2) $\alpha_{u^{(l)},s}(I)\le M^{-(s-1)/2}$ if 
$1\le s\le d-|I|-1$. This means that formula~(5.6) holds for the 
operator we get by an appropriate reindexation the indices 
$\{1,\dots,d\}\setminus I$ of the arguments of $B^I_{u^{(l)}}$ to 
the set $1,\dots,d-|I|$. An appropriate reindexation is obtained 
if the elements of the set $\{1,\dots,d\}\setminus I$ are listed 
in a monotone increasing order, and the $j$-th element of this 
sequence gets the index~$j$. Such a reindexation of the indices 
yields a version of Lemma~6.3 that enables us to estimate the 
terms $W_K^u(\alpha_{u^{(l)}}^I,t)$. (Originally we get an estimate
for a version of $W_K^u(\alpha_{u^{(l)}}^I,t)$ with reindexed
parameters by means  of a version of $B^I_{u^{(l)}}$ with reindexed
parameters.) 

In the application of Lemma~6.3 we still have to understand what
$\tilde\alpha_{j,k}(u_j)$ means in formula~(6.9) if $B^I_{u^{(l)}}$ 
plays the role of the operator~$A$.
 
By formula (6.8) in Lemma~6.3 we get that
$$
W^u_K(\alpha^I_{u^{(l)}},t)\le\frac{t^{|K|}}
{M^{(d-|I|-|K|-1)/2}}\le c^2 2^{-2N}M^{-(d-|I|-1)/2}
\quad \text{if }2\le|K|\le d-|I|-1.
$$
I claim that relations (6.9) and (7.3) imply that
$$
W^u_{\{k\}}(\alpha^I_{u^{(l)}},t)\le t
\tilde\alpha^I_{k}(u^{(l)})\le c2^{-2N}M^{-(d-|I|-1)/2}
$$
for a one point set $\{k\}\in\{1,\dots,d-1\}\setminus I$. We get 
this bound from~(7.3) if we show that
$\tilde\alpha_{j,k}(u_j)\le \alpha^I_k(u^{(l)})$ for any 
$j\in\{1,\dots,d-1\}\setminus I$, $j\neq k$, with the function 
$\tilde\alpha_{j,k}(u_j)$ corresponding to the operator 
$B^I_{u^{(l)}}$ if $u_j\in B^{n_j}$.
 
This inequality can be seen by giving a good representation
of $\tilde\alpha_{j,k}(u_j)$ 
when it corresponds to $B^I_{u^{(l)}}$
instead of $A$ together with a similar representation of
$\tilde\alpha^I_k(u^{(l)})$. An adaptation of formula~(5.4) will
be applied to this case. The main difference between formula~(5.4)
and the representation of $\tilde\alpha_{j,k}(u_j)$ given below is
that in the new formula we have the fixed functions
$u^{(l)}_s(\cdot)$ in the coordinates $s\in I$. In this case we have
$$
\align
\tilde\alpha_{j,k}(u_j)&=\sup_{v_p(\cdot),\,p\in\{1,\dots,d-1\}
\setminus(I\cup\{j,k\},\;v_{k,d}(\cdot,\cdot)}
\sum_{i_1,\dots,i_d} a(i_1,\dots,i_d) u_j(i_j)v_{k,d}(i_k,i_d)\\
& \qquad \prod_{s\in I} u_s^{(l)}(i_s)
\prod_{p\in\{1,\dots,d-1\}\setminus (I\cup\{j,k\})} v_p(i_p)
\endalign
$$
for a vector $u_j\in R^{n_j}$, where the supremum is taken for
such vectors $v_p(\cdot)$ depending on
the coordinate $i_p$, $p\in\{1,\dots,d-1\}\setminus(I\cup\{j,k\})$,
for which $\summ_{i_s}v_p^2(i_p)\le1$ and a function
$v_{k,d}(\cdot,\cdot)$, depending on the coordinates $i_k$ and
$i_d$ such that $\summ_{i_k,i_d}v^2(i_k,i_d)\le1$. The expression
$\tilde\alpha^I_k(u^{(l)})$ has a similar representation,
only in its definition we have to take supremum also for all
vectors $v_j(\cdot)\in B^{n_j}$ in its $j$-th coordinate
instead of fixing a vector $u_j\in B^{n_j}$ as it was done in the
definition of $\tilde\alpha_{j,k}(u_j)$, $u_j\in B^{n_j}$. These
observations imply the desired inequality
$\tilde\alpha_{j,k}(u_j)\le \alpha^I_k(u^{(l)})$.
 
The above inequalities imply that
$$
\sum_{J\colon\;J\subset\{1,\dots,d-1\}\setminus I,\,J\neq\emptyset}
W^u_J(\alpha^I_{u^{(l)}},4t)\le 2^{-2N}M^{-(d-|I|-1)/2}
$$
for all $u\in B^{n_{j_1}}\times\cdots\times B^{n_{j_s}}$
if the constant $c>0$ in the choice $t=c2^{-N}M^{-1/2}$ is
sufficiently small. Hence it follows from the corollary of
Proposition~4.2 applied for the metric
$\bar\rho_{\alpha^I_{u^{(l)}}}$ induced by the norm 
$\alpha^I_{u^{l)}}$ with the choice
$D=B^{n_{j_1}}\times\cdots\times B^{n_{j_r}}$ and $t=c2^{-N}M^{-1/2}$
with a sufficiently small number $c>0$ and 
$u=2^{-2N}M^{-(d-|I|-1)/2}$ that relation~(7.1) holds.
 
In the case $|I|=d-2$ we can write $I=\{1,\dots,d-1\}\setminus\{k\}$
with an appropriate $k\in\{1,\dots,d-1\}$. The inequality
$\tilde\alpha_{j,k}(u^{(l)})\le 2^{-N}M^{(d-2)/2}$ with an 
arbitrary index $j\in I$ implies in this case that
$$
\sum_{i_k,i_d} b^I_{u^{(l)}}(i_k,i_d)v(i_k,i_d)
\le M^{|I|/2}2^{-N}M^{-(d-2)/2}\le 2^{-N}
\quad \text{if }\summ_{i_k,i_d}v^2(i_k,i_d)\le1,
$$
or in an equivalent form
$$
\summ_{i_k,i_d} b^I_{u^{(l)}}(i_k,i_d)^2\le 2^{-2N}, \tag7.4
$$
where $I=\{1,\dots,d-1\}\setminus\{k\}$, and the numbers 
$b^I_{u^{(l)}}(i_k,i_d)$ are defined in~(6.1). Let us also define 
the pseudonorm
$$
\beta^I_{u^{(l)}}(v)=
\[\sum_{i_d}\(\sum_{i_k} b^I_{u^{(l)}}(i_k,i_d)v(i_k)\)^2\]^{1/2}
$$
of the vectors $v=(v(1),\dots,v(n_k))\in R^{n_k}$.
 
The pseudometric $\bar\rho_{\alpha^I_{u^{(l)}}}$ defined in (6.5) 
agrees in this case with the metric induced by the pseudonorm
$\beta^I_{u^{(l)}}$. Hence in this case the existence of a partition
$V_1,\dots,V_L$  of $B^{n_k}$ with $L\le 2^{C2^{2N}M}$ elements
and the property
$\bar\rho_{\alpha_{u^{(l)}}^I}(v,\bar v)\le2^{-2N}M^{-1/2}$,
if $v,\bar v\in V_l$ with some $1\le l\le L$, i.e. relation~(7.1) 
can be proved with the help of the corollary of Proposition~4.2 
and the following estimate on the pseudonorm $\beta^I_{u^{(l)}}$.
 
By the Schwarz inequality and formula (7.4)
$$
\aligned
E\beta^I_{u^{(l)}}(G_k)&\le
\[\sum_{i_d}E\(\sum_{i_k} b^I_{u^{(l)}}(i_k,i_d)g_k(i_k)\)^2\]^{1/2}\\
&=\[\sum_{i_k,i_d}b^I_{u^{(l)}}(i_k,i_d)^2\]^{1/2}\le 2^{-N}
\endaligned \tag7.5
$$
for a standard normal random vector $G_k=(g_k(1),\dots,g_k(n_k))$
of dimension $n_k$. Because of relation~(7.5) an application of the 
corollary of Proposition~4.2 for the operator $\beta^I_{u^{(l)}}$ 
with $4t=2^{-(N+1)}M^{-1/2}$ and $u=2^{-2N}M^{-1/2}$ shows
the existence of a partition $V_1,\dots,V_L$ of $B^{n_k}$ with 
$L\le 2^{C_1/t^2}\le 2^{C2^{2N}M}$ elements such that
$\bar\rho_{\alpha_{u^{(l)}}^I}(v,\bar v)=\beta^I_{u^{(l)}}(v-\bar v)
\le2^{-2N}M^{-1/2}$ if $v\in V_l$, $\bar v\in V_l$ with some
$1\le l\le L$. We had to prove this statement.
 
Let us fix some $u^{(l)}$ appearing in the partition $u^{(l)}+U_l$,
$1\le l\le L$ of the set $U$ we are considering. It can be shown
with the help of relation (7.1) that there exists a partition
$V_1(l),\dots,V_{L_l}(l)$ of $B^{n_1}\times\dots\times B^{n_{d-1}}$
with $L_l\le 2^{C2^{2N}M}$ elements such that
$$
\aligned
\rho_{\alpha_{u^{(l)}}^I} (v,\bar v)&\le2^{-2N}M^{-(d-|I|-1)/2}
\quad \text{if } v\in V_p(l) \text{ and } v\in V_p(l)
\text{ with the same} \\ 
&\qquad \text{index } p \text{ for all } I\subset\{1,\dots,d-1\}
\text{ such that } 1\le |I|\le d-2.
\endaligned\tag7.6
$$
 
Indeed, it follows from (7.1) and the definition of
$\rho_{\alpha_{u^{(l)}}^I}$ in~(6.5) and~(6.7) that for all sets
$I\subset\{1,\dots,d-1\}$, $1\le |I|\le d-2$, there is a partition
$V_1(l,I)$,\dots, $V_{L(l,I)}(l,I)$ of
$B^{n_1}\times\dots\times B^{n_{d-1}}$ depending on $l$ and $I$
with $L(l,I)\le 2^{C2^{2N}M}$ elements such that
$\rho_{\alpha_{u^{(l)}}^I} (v,\bar v)\le2^{-2N}M^{-(d-|I|-1)/2}$
if $v\in V_p(l,L)$ and $\bar v\in V_p(l,L$ with the same index~$p$.
Then taking all possible intersections
$\bigcapp_{I\colon\;I\in\{1,\dots,d-1\},\,1\le|I|\le d-2}V_{p(I)}(I,l)$
that contain exactly one element from each above introduced
partitions depending on the sets $I\subset\{1,\dots,d-1\}$, 
$1\le |I|\le d-2$, we get a partition of 
$B^{n_1}\times\dots\times B^{n_{d-1}}$ that satisfies (7.6). Let us 
observe that the number of elements of this partition also can be 
bounded from above by $2^{C2^{2N}M}$ with some constant $C>0$.
 
Let us choose a partition $V_1(l),\dots,V_{L_l}(l)$ of
$B^{n_1}\times\dots\times B^{n_{d-1}}$ satisfying relation~(7.6)
for all vectors $u^{(l)}$ taking part in a partition
$u^{(l)}+U_l$, $1\le l\le L$ satisfying Lemma~6.1. Then the ensemble
of sets $u^{(l,p)}+(V_p(l)\cap U_l)$, $1\le p\le L(l)$,
$1\le l\le L$, with $u^{(l,p)}=u^{(l)}$ constitutes a partition
of the set $U$ which, after an appropriate reindexation,
satisfies Lemma~6.2.
 
\beginsection 8. The proof of the Main inequality.
 
In this section I prove the Main inequality with the help of Lemma~6.2.
 
\medskip\noindent
{\it The proof of the Main inequality.}\/ First it will be shown that
relation (5.7) holds with an appropriate constant $C=C(d)$ in it
if $N\ge N_0$ with a sufficiently large threshold index $N_0=N_0(r)$.
To this end let us observe that
$$
E(Y_d(u)-Y_d(u'))^2=\rho_\alpha(u,u')^2,
$$
hence
$$
E(Y_d(u)-Y_d(u'))^{2M}=1\cdot3\cdot\cdots\cdot(2M-1)\rho_\alpha(u,u')^{2M}
\le(2M)^M\rho_\alpha( u,u')^{2M}
$$
with the metric $\rho_\alpha$ defined in (2.11) for arbitrary vectors
$u\in R^{n_1}\times\cdots\times R^{n_{d-1}}$,
$u'\in R^{n_1}\times\cdots\times R^{n_{d-1}}$ and $M\ge1$.
In particular,
$$
\aligned
E(Y_d(u^{(t)})-Y_d(u^{(t')}))^{2^{2(N+A)}M}
&\le(2^{2(N+A)}M)^{2^{2(N+A)}M/2}
\cdot (2^{-2N}M^{-(d-1)/2})^{2^{2(N+A)}M}\\
&=(2^{(A-N)}M^{-(d-2)/2})^{2^{2(N+A)}M}
\endaligned
\tag8.1
$$
for all $u^{(t)}\in U$ and $u^{(t')}\in U$ if
$U\in\Cal U(r,N)$. As a consequence,
$$
\align
E&\[\sup_{(u^{(t)},u^{(t')})\colon\; u^{(t)}\in U,\,
u^{(t')}\in U}(Y_d(u^{(t)})-Y_d(u^{(t')}))\]^{2^{2(N+A)}M}\\
&\qquad \le r^2 (M^{-(d-2)/2}2^{(A-N)})^{2^{2(N+A)}M}
\le(2M^{-(d-2)/2}\cdot 2^{(A-N)})^{2^{2(N+A)}M}\\
&\qquad\le (CM^{-(d-2)/2}2^{(A-N)})^{2^{2(N+A)}M}
\endalign
$$
if $N\ge N_0(r)$ with some threshold $N_0(r)$ and constant
$C\ge2$, i.e. relation (5.7) holds for $N\ge N_0$ with
$C=C(d)\ge2$ and $A\ge A_0\ge0$.
 
Hence it is enough to show that relation (5.7) holds for a
set $U\in\Cal U(r,N)$ if it holds for all sets $U\in\Cal U(r,N+2)$.
To show this let us consider such a partition $u^{(l)}+U_l$,
$1\le l\le L$ of the set $U\in\Cal U(r,N)$ with
$L\le 2^{C2^{2N}M}$ elements which satisfies Lemma~6.2. First the 
following weaker estimate will be verified.
 
Let us take an element $u^{(l)}+U_l$ of the partition of $U$ we
consider. Let us denote this set by $\bar U_l$. We will show that
the estimate
$$
\aligned
&E\[\sup_{(u^{(t)},u^{(t')})\colon\; u^{(t)}\in \bar U_l,\,
u^{(t')}\in \bar U_l}(Y_d(u^{(t)})-Y_d(u^{(t')}))\]^{2^{2(N+A)}M} \\
&\qquad\le \(\frac C3M^{-(d-2)/2}2^{(A-N)}\)^{2^{2(N+A)}M}
\endaligned \tag8.2
$$
holds for all $A\ge A_0$ with some threshold index $A_0$ and the
same constant $C=C(d)$ which appears in (5.7) (with parameter $N+2$)
if these constant (depending only on the parameter~$d$) are chosen
sufficiently large.
 
To prove relation (8.2) let us consider two arbitrary vectors
$u\in\bar U_l$ and $u'\in\bar U_l$, write them in the
form $u=u^{(l)}+u^{(0)}$ and $u'=u^{(l)}+{u'}^{(0)}$
with $u^{(0)}\in U_l$ and ${u'}^{(0)}\in U_l$. We can write the
difference $Y_d(u)-Y(u')$ because of the special form of relation~(2.1)
defining $Y_d(u)$ as
$$
\aligned
Y_d(u)-Y_d(u')&=Y_d(u^{(l)}+u^{(0)})-Y_d(u^{(l)}+ {u'}^{(0)})
=Y_d(u^{(0)})-Y_d({u'}^{(0)})  \\
&\qquad +\sum_{I\colon\;I\subset\{1,\dots,d-1\},\,1\le |I|\le d-2}
M^{-|I|/2}\[Y^I_{u^{(l)}}(u^{(0)})-Y^I_{u^{(l)}}({u'}^{(0)})\],
\endaligned \tag8.3
$$
where
$$
Y^I_{u^{(l)}}(v)=\sum_{(i_j,\,j\in \{1,\dots,d-1\}\setminus I}
b^I_{u^{(l)}}(i_j,\,j\in\{1,\dots,d\}\setminus I)
\prod_{j\in\{1,\dots,d-1\}\setminus I} v_j(i_j)g_d(i_d)
$$
for all
$v=(v_1(i_1),\dots,v_{d-1}(i_{d-1}))\in R^{n_1}\times\cdots\times R^{n_{d-1}}$
and $I\subset\{1,\dots,d-1\}$, $1\le |I|\le d-2$ with the constants
$b^I_{u^{(l)}}(i_j,\,j\in\{1,\dots,d\}\setminus I)$ defined
in (6.1). (Here we apply this formula with the choice $u=u^{(l)}$,)
and $(g_d(1),\dots,g_d(n_d))$ is the same vector of independent,
standard Gaussian random variables which appeared in the definition
of $Y_d(u)$.)
 
In the subsequent considerations the following notation will be applied.
Given some vector $u^{(t)}\in \bar U_l$, its decomposition to the vector
$u^{(l)}$ plus a vector in $U_l$ will be denoted as
$u^{(t)}=u^{(l)}+u^{(t,0)}$ with $u^{(t,0)}\in U_l$.
 
By taking the supremum of the expressions both at the left-hand
and right-hand side of identity~(8.3) for all pairs
$(u^{(t)},u^{(t')})$ such that $u^{(t)}\in\bar U_l$ and
$u^{(t')}\in\bar U_l$ we get an identity that implies the following
inequality.
$$
\sup_{(u^{(t)},u^{(t')})\colon\; u^{(t)}\in \bar U_l,\,
u^{(t')}\in \bar U_l}(Y_d(u^{(t)})-Y_d(u^{(t')}))\le Z+
\sum_{I\colon\; I\subset \{1,\dots,d-1\}, \;1\le|I|\le d-2}
M^{-|I|/2} Z_I \tag8.4
$$
with
$$
Z=Z(l,N)=\sup_{(u^{(t,0)},u^{(t',0)})\colon\; u^{(t)}\in U_l,\,
u^{(t',0)}\in U_l}(Y_d(u^{(t,0)})-Y_d(u^{(t',0)}))
$$
and
$$
\aligned
Z_I=Z_I(l,N)&=\sup_{(u^{(t,0)},u^{(t',0)})\colon\; u^{(t)}\in U_l,\,
u^{(t',0)}\in U_l}
[Y^I_{u^{(l)}}(u^{(t,0)})-Y^I_{u^{(l)}}(u^{(t',0)})], \\
&\qquad \text{for all } I\subset\{1,\dots,d-1\}
\text{ such that } 1\le |I|\le d-2.
\endaligned
$$
I claim that
$$
EZ^{2^{2(N+A)}M}\le\(\frac C4M^{-(d-2)/2}2^{(A-N)}\)^{2^{2(N+A)}M} \tag8.5
$$
with the same constant $C=C(d)$ as in formula (5.7) if $A\ge A_0$ with
some fixed number $A_0=A_0(d)\ge0$, and
$$
EZ_I^{2^{2(N+A)}M}\le(C'M^{-(d-|I|-2)/2}2^{(A-N)})^{2^{2(N+A)}M}
\tag8.6
$$
for all $I\subset\{1,\dots,d-1\}$ such that $1\le |I|\le d-2$
if $A\ge A_0$ with some universal constants $A_0$ and $C'$.
Let me emphasize in particular that the constants $A_0$ and $C'$
in~(8.6) do not depend on the choice of the constant $C=C(d)$ and 
threshold index~$A_0$ in~(5.7).

Relation (8.5) can be deduced from our inductive hypothesis by
which the Main inequality holds for $N+2$ and the fact that
$U_l\in\Cal U(r,N+2)$. Indeed, this inductive hypothesis together
with H\"older's inequality yield that
$$
\aligned
&EZ^{2^{2(N+A)}M}=E\[\sup_{((u^{(t,0)},u^{(t',0)}))\colon\;
u^{(t,0)}\in U_l,\,u^{(t',0)}\in U_l}
(Y_d((u^{(t,0)})-Y_d(u^{(t',0)})\]^{2^{2(N+A)}M}\\
&\qquad\le\[ E\[\sup_{((u^{(t,0)},u^{(t',0)}))\colon\;
u^{(t,0)}\in U_l,\,u^{(t',0)}\in U_l}
(Y_d((u^{(t,0)})-Y_d(u^{(t',0)})\]^{2^{2(N+A+2)}M}\]^{1/4}\\
&\qquad\le\(CM^{-(d-2)/2}2^{(A-(N+2))}\)^{2^{2(N+A+2)M}/4}
=\(\frac14CM^{-(d-2)/2}2^{(A-N)}\)^{2^{2(N+A)}M}.
\endaligned
$$
(The reason for applying the induction from $N+2$ and not from $N+1$ to
$N$ in our proof is that in such a way we got a coefficient $\frac14$
at the right-hand side of estimate~(8.5). An induction from $N+1$ to
$N$ would yield only a weaker estimate with multiplying factor $\frac12$
which would be not sufficient for our purposes.)
 
Relation (8.6) will be proved first only in the case $1\le |I|\le d-3$. 
This will be done with the help of the Main inequality with parameter 
$d-|I|\le d-1$. This is legitime because of our inductive hypothesis. 
The main inequality will be applied for the operator $B^I_{u^{(l)}}$ 
defined in~(6.2) as the operator $A$ and the set of vectors 
$U\in\Cal U(r,N)$ will be chosen as 
$U=U_l(I)=\{u^{(t,0)}_{I^c}\colon u^{(t,0)}\in U_l\}$. That is we get 
the set $U$ by taking the vectors $u=(u_1,\dots,u_{d-1})\in U_l$
and omitting their coordinates indexed by the elements of the
set~$I$. More precisely, we apply the Main inequality for a version
of $B^I_{u^{(l)}}$ and $U_l(I)$ we get by renumerating the indices
of their coordinates to the sets $\{1,\dots,d-|I|\}$ and
$\{1,\dots,d-|I|-1\}$ respectively in an appropriate way. A good
way of reindexation of the coordinates is to list them with monotone
increasing indices and to give then the $j$-th element the
index~$j$.

To apply the Main inequality we have to show that its conditions 
are satisfied with such a choice. We have to check that the operator 
$B^I_{u^{(l)}}$ satisfies relation~(5.6). (Here $d-|I|$ takes the 
role of the parameter~$d$.) This statement follows from the analogous 
statement for the operator~$A$. Beside this, we have to show that
$U_l(I)\in\Cal U_{B^I_{u^{(l)}},M,d-|I|}(r,N)$. This can be
done with the help of Lemma~6.2.

The estimate we have to give on 
$\bar\rho_{\alpha^I_{u^{(l)}}}(\cdot,\cdot)$ to show that
$U_l(I)\in\Cal U(r,N)=\Cal U_{B^I_{u^{(l)}},M,d-|I|}(r,N)$
agrees with the estimate we proved in Lemma~6.2 on this quantity.
The bound we have to give about $\tilde\alpha_{j,k}$  to show that
$U_l(I)\in\Cal U_{B^I_{u^{(l)}},M,d-|I|}(r,N)$ follows from
relation~(7.3) and the inequality
$\tilde\alpha_{j,k}(u_j)\le\tilde\alpha_k^I(u^{(l)})$
if $u_j\in B^{n_j}$ with the same quantities 
$\tilde\alpha_{j,k}(u_j)$ and $\tilde\alpha_k^I(u^{(l)})$ which 
appeared in the proof of Lemma~6.2. The remaining properties
needed to check that $U_l(I)\in\Cal U_{B^I_{u^{(l)}},M,d-|I|}(r,N)$
clearly hold. Then the Main inequality may be applied with such 
a choice, and it yields relation~(8.6) for $1\le |I|\le d-3$.
 
In the case $|I|=d-2$, and the set $\{1,\dots,d-1\}\setminus I$
consists of a point $k$, and formula (8.6) can be
proved with the help of the Main inequality in the case $d=2$ in
a similar way. This inequality can be applied for the operator
$2^NB^I_{u^{(l)}}$ defined as $2^NB^I_{u^{(l)}}(v)
=\sum\limits_{i_k,i_d}2^{N} b^I_{u^{(l)}}(i_k,i_d)v(i_k,i_d)$
for a vector $v\in R^{n_k}\otimes R^{n_d}$ as the operator~$A$.
It follows from Lemma~6.2 that
$$
\alpha_2(u^{(t,0)}-u^{(t',0)})=
\[\sum_{i_d}\(\sum_{i_k}2^N b^I_{u^{(l)}}(i_k,i_d)
[u^{(t,0)}_{I^c}(i_k)-u^{(t',0)}_{I^c}(i_k)]\)^2\]^{1/2}
\!\!\!\le2^{-N}M^{-1/2}
$$
if $u^{(t,0)}(i_k)\in U_l$ and $u^{(t',0)}(i_k)\in U_l$. Hence the
set consisting of all vectors of the form
$\frac12(u^{(t,0)}_{I^c}-u^{(t,0)}_{I^c})$ with some
$u^{(t,0)}\in U_l$ and $u^{(t',0)}\in U_l$ is contained in the
set $U_N$ introduced in~(3.2). The
inequality $\alpha_1(2^NB_{u^{(l)}}^I)\le1$ also holds by formula~(7.4).
Hence the Main inequality in the case $d=2$ can be applied in
this case, and it yields that
$E(2^NZ_I)^{2^{2(N+A)}M}\le(C\cdot2^A)^{2^{2(N+A)}M}$
which is equivalent to (8.6) for $|I|=d-2$.
 
Inequality (8.2) follows from relations (8.4), (8.5), (8.6) and
Minkowski's inequality for $L_p$ norms with $p=2^{2(N+A)}M$. (Observe
that we are working with non-negative random variables, since the
supremums we consider contain the terms
$Y_d(u^{(t)})-Y_d(u^{(t)})\equiv0$.) Indeed, they yield that
$$
\aligned
&E\[\sup_{(u^{(t)},u^{(t')})\colon\; u^{(t)}\in \bar U_l,\,
u^{(t')}\in \bar U_l}(Y_d(u^{(t)})-Y_d(u^{(t')}))\]^{2^{2(N+A)}M} \\
&\qquad\le \(\(\frac C4+2^dC'\)M^{-(d-2)/2}2^{(A-N)}\)^{2^{2(N+A)}M}.
\endaligned
$$
If the constant $C=C(d)$ in the Main inequality is chosen sufficiently
large, then $\frac C4+2^dC'\le\frac C3$ in the last inequality, and
this means that it implies relation~(8.2).
 
The Main inequality will be proved with the help of inequality~(8.2).
It will be also exploited that the cardinality of the partition of a 
set~$U$ in Lemma~6.2 is not too large.
 
Let us consider a partition $\bar U_l=u^{(l)}+U_l$, $1\le l\le L$,
of a set $U\in\Cal U(r,N)$ with $L\le 2^{C2^{2N}M}$ elements that
satisfies Lemma~6.2. Let us fix an element $\bar u^{(l)}\in\bar U_l$
in all sets $\bar U_l$, $1\le l\le L$. Given a vector
$u^{(t)}\in U$ let $\ell(t)$ denote that index $l$, $1\le l\le L$,
for which $u^{(t)}\in\bar U_l$. Then we can write for two arbitrary
vectors $u^{(t)}\in U$ and $u^{(t')}\in U$ the inequality
$$ \allowdisplaybreaks
\align
&Y_d(u^{(t)})-Y_d(u^{(t')})=[Y_d(u^{(t)})-Y_d(\bar u^{\ell(t)})]+
[Y_d(u^{(t')})-Y_d(\bar u^{\ell(t')})]\\
&\qquad\qquad +[Y_d(\bar u^{\ell(t')})-Y_d(\bar u^{\ell(t)})] \\
&\qquad \le2\sup_{1\le l\le L} \sup_{u^{(s)}\colon\;u^{(s)}\in\bar U_l}
[Y_d(u^{(s)})-Y_d(\bar u^{(l)}]
+\sup_{1\le l,l'\le L} [Y_d(\bar u^{(l')})-Y_d(\bar u^{(l)})].
\endalign
$$
Since the right-hand side of the above inequality does not depend
on the vectors $u^{(t)}\in U$ and $u^{(t')}\in U$ it implies that
$$
\aligned
&\sup_{(u^{(t)},u^{(t')})\colon\;u^{(t)}\in U,\; u^{(t')}\in U}
[Y_d(u^{(t)})-Y_d(u^{(t')})]\\
&\qquad \le2\sup_{1\le l\le L} \sup_{u^{(s)}\colon\;u^{(s)}\in\bar U_l}
[Y_d(u^{(s)})-Y_d(\bar u^{(l)}]
+\sup_{1\le l,l'\le L} [Y_d(\bar u^{(l')})-Y_d(\bar u^{(l)})].
\endaligned \tag8.7
$$
 
The Main inequality can be proved by means of good moment estimates
on the two terms at the right-hand side of inequality~(8.7).
It follows from inequalities (8.2) and $L\le2^{C'2^{2N}M}$ for
the number of partitions of~$U$ in Lemma~6.2, where the number
$C'$ does not depend on the constants $A_0=A_0(d)$ and $C=C(d)$
in the Main inequality that
$$
\aligned
&E\(2\sup_{1\le l\le L} \sup_{u^{(s)}\colon\;u^{(s)}\in\bar U_l}
[Y_d(u^{(s)})-Y_d(\bar u^{(l)}]\)^{2^{2(N+A)}M} \\
&\qquad \le\sum_{l=1}^L E\(2\sup_{(u^{(s)},u^{(s')})
\colon\;u^{(s)}\in\bar U_l,\; u^{(s')}\in U_l}
[Y_d(u^{(s)})-Y_d(\bar u^{(s')}]\)^{2^{2(N+A)}M} \\
&\qquad\le L\(\frac{2C}3M^{-(d-2)/2}2^{(A-N)}\)^{2^{2(N+A)}M} \\
&\qquad \le 2^{C'2^{2N}M}\(\frac{2C}3M^{-(d-2)/2}
2^{(A-N)}\)^{2^{2(N+A)}M} \\
&\qquad =\(2^{C'2^{-2A}}\frac{2C}3M^{-(d-2)/2}2^{(A-N)}\)^{2^{2(N+A)}M}
\le\(\frac{3C}4M^{-(d-2)/2}2^{(A-N)}\)^{2^{2(N+A)}M}
\endaligned \tag8.8
$$
if the threshold index $A_0$ in the Main inequality is chosen so large
that $2^{C'2^{-2A}}\le\frac98$ for $A\ge A_0$. Such a choice is
possible since the constant $C'$ appearing in the exponent of the
bound for the cardinality of the partition of the set $U$ does not
depend on the choice of the number~$A_0$ in the Main inequality.
(The threshold index $A_0$ was introduced to have a control on the
multiplicative factor $L$ in the previous estimate.)
 
We get in a similar way with the help of inequality (8.1)
$$ \allowdisplaybreaks
\align
&E\(\sup_{1\le l,l'\le L}
[Y_d(\bar u^{(l')})-Y_d(\bar u^{(l)})]\)^{2^{2(N+A)}M}
\le\sum_{1\le l,l'\le L}
E\([Y_d(\bar u^{(l')})-Y_d(\bar u^{(l)})]\)^{2^{2(N+A)}M}\\
&\qquad \le L^2\(M^{-(d-2)/2}2^{(A-N)}\)^{2^{2(N+A)}M}
\le 2^{2C'2^{2N}M}\(M^{-(d-2)/2}2^{(A-N)}\)^{2^{2(N+A)}M}\\
&\qquad \le \(\bar C M^{-(d-2)/2}2^{(A-N)}\)^{2^{2(N+A)}M}
\tag8.9
\endalign
$$
with a constant $\bar C$ which does not depend on the constant
$C=C(d)$ in the Main inequality.
 
It follows from relations (8.7), (8.8), (8.9) and Minkowski's
inequality for $L_p$~norms with $p=2^{2(N+A)}M$ (we are working
again with non-negative random variables) that
$$
\align
&E\[\sup_{(u^{(t)},u^{(t')})\colon\; u^{(t)}\in U,\;u^{(t')}\in U}
(Y_d(u^{(t)})-Y_d(u^{(t')})) \]^{2^{2(N+A)}M} \\
&\qquad \le\( \(\frac{3C}4+\bar C\)M^{-(d-2)/2}2^{(A-N)}\)^{2^{2(N+A)}M}
\le(CM^{-(d-2)/2}2^{(A-N)})^{2^{2(N+A)}M}
\endalign
$$
if the constant $C=C(d)$ (together with $A_0=A_0(d)$) is chosen
sufficiently large. The Main inequality is proved.
 
\medskip\noindent
{\it Remark.}\/ In the proof of the Main inequality Lemma~6.2 played
an important role. In the verification of relation~(8.2) we have
exploited that the partition of the set $U$ we have considered
satisfies all properties formulated in Lemma~6.2. Lata{\l}a tried
to prove an analogous estimate for all partitions satisfying
Lemma~6.1. His proof however contains an error. It applies an
estimate formulated in relation~(18) of Lemma~7 in paper~[2]
whose verification is based on Theorem~3 of~[2]. But the proof of
this Theorem~3 is incorrect. This result was proved by means of
a backward induction similarly to our Main inequality. In
Lata{\l}a's backward induction procedure we turn from the
parameters $(r-1,l+1)$ to $(r,l)$. (In this remark I apply the
notation of paper~[2].) But the last step of this backward induction
when we turn from $l=1$ to $l=0$ does not work.
 
The essential statement of Theorem~3 in~[2] proved by backward
induction was formulated by means of some objects denoted by
$\Delta_l$ and $\tilde\Delta_l$. These objects were defined
differently for $l=0$ and $l\ge1$. Hence a special argument would
have been needed to prove the induction step (formulated in
relation~(25) of~[2]) for $l=0$. But such a step is missing from
the proof of~[2]. Moreover, it would require different arguments.
The situation seems to be similar to the proof of the Basic estimate
with the help of the Main inequality in this paper
where formula~(5.6) had to be replaced by the stronger condition~(1.10).
I think that in the last step of the backward induction proof of
Theorem~3 in~[2] such a bound should appear which also depends on
the term $\|A\|_{\{1\},\dots,\{d\}}$, and this would supply only a
weaker estimate.
 
This seems to be a serious error. I believe that
not only the proof of Theorem~3 is incorrect, but even the results
formulated in Theorem~3 and relation~(18) of~[2] are wrong.
Since Lata{\l}a's proof heavily exploited formula~(18) it was not
clear for me whether his main result holds in its original form or
it must be modified. My main goal in this paper was to answer this
question. Finally it turned out that Lata{\l}a's result is correct.
But to prove this I had to find a new, better partition of the sets
$U\in\Cal U(r,N)$ than Lata{\l}a did. It was the partition
constructed in Lemma~6.2 that helped in saving Lata{\l}a's proof.

\beginsection Appendix. The proof of Propositions 4.1 and 4.2.

{\it The proof of Proposition~4.1.}\/ Put 
$K=\{y\colon\;y\in R^n,\; \alpha_1(y)\le4\alpha_1(tG),\; 
\alpha_2(y)\le4\alpha_2(tG)\}$. Then $\mu_{n,t}(K)\ge\frac12$, 
since by the Markov inequality 
$$
1-\mu_{n,t}(K)\le\mu_{n,t}(y\colon\;\alpha_1(y)>4tE\alpha_1(tG))
+\mu_{n,t}(y\colon\;\alpha_2(y)>4tE\alpha_2(tG))\le\frac12.
$$ 
Beside this, the set $K$ has the symmetry property $-K=K$ which 
yields that
$$
\align
&\mu_{n,t}(y\colon\; y\in R^n,\;\alpha_1(y-x)\le 4E\alpha_1(tG),\;
\alpha_2(y-x)\le 4E\alpha_2(tG))\\
&\qquad=C_n\int_{K+x}e^{-y^2/2t}\,dy=C_n\int_Ke^{-(y+x)^2/2t}\,dy
=e^{-\|x\|^2/2t}\int_Ke^{(y,x)/t}\mu_{n,t}(\,dy)\\
&\qquad=e^{-\|x\|^2/2t}\int_K\frac12\(e^{(y,x)/t}+e^{(-y,x)/t}\)
\mu_{n,t}(\,dy)\ge e^{-\|x\|^2/2t}\mu_{n,t}(K)
\endalign
$$
with the norming constant $C_n=(\sqrt{2\pi} t)^{-n}$. Hence the relations
$\mu_{n,t}(K)\ge\frac12$, and $\|x\|\le1$ (i.e. $x\in B^n$) imply that
$$
\align
&\mu_{n,t}(\{y\colon\; y\in R^n,\;\alpha_1(y-x)\le 4E\alpha_1(tG),\;
\alpha_2(y-x)\le 4E\alpha_2(tG)\})\\
&\qquad \ge \frac12e^{-\|x\|^2/2t}\ge \frac12e^{-1/2t^2}.
\endalign
$$

\medskip\noindent
{\it The proof of Proposition 4.2.}\/ In the case $d=1$
Proposition~4.2 immediately follows from Proposition~4.1 if it is
applied for $\alpha=\alpha_1=\alpha_2$, and the relation
$4E\alpha(tG_n)=E\alpha(4tG_n)$ is exploited. Hence it is enough
to prove Proposition~4.2 for $d$ under the inductive hypothesis
that it holds for $d-1$.
 
Let us fix some
$x=(x_1,\dots,x_d)\in B^{n_1}\times\cdots\times B^{n_d}$, where $B^n$
denotes the unit ball in~$R^n$. We can write
$$
\aligned
\rho_\alpha(x,y)&=\alpha(y_1\otimes\cdots\otimes y_{d-1}\otimes y_d
-x_1\otimes\cdots\otimes x_{d-1}\otimes x_d)\\
&\le \alpha(x_1\otimes\cdots\otimes x_{d-1}\otimes(y_d-x_d))+
\alpha((y_1\otimes\cdots\otimes y_{d-1}-x_1\otimes\cdots x_{d-1})\otimes y_d)
\endaligned \tag A1
$$
for arbitrary $y=(y_1,\dots,y_d)\in R^{n_1}\times\cdots\times R^{n_d}$.
 
We shall define some sets $A$, $B$ and $C$. We shall not denote their
dependence on the vector
$x\in B^{n_1}\times\cdots\times B^{n_d}$ we have fixed. To define
the set $A$ first we introduce the following
quantity $W_I^x(y|\alpha,t)$ similar to the quantity
$W_I^x(\alpha,t)$ defined in formula~(4.2).
 
Let us fix $d$ independent standard normal vectors
$G_j=(g_j(1),\dots,g_j(n_j))$ of dimension $n_j$, $1\le j\le d$,
and define for all $t>0$, $y\in R^{n_d}$ and $I\subset\{1,\dots,d\}$,
$I\neq\emptyset$ the quantity
$$
\align
W_I^x(y|\alpha,t)&=E\alpha(z_1\otimes\cdots\otimes z_d)\quad
\text{where } z_j=tG_j \text{ if } j\in I, \\
&\qquad z_j=x_j \text{ if }j\notin I \text{ and } j\neq d,
\quad \text{and } z_d=y \text{ if } d\notin I.
\endalign
$$
Then we put
$$
\align
A&=\biggl\{y_d\colon\; y_d\in R^{n_d},
\;\alpha(x_1\otimes\cdots\otimes x_{d-1}\otimes(y_d-x_d))
\le E\alpha(x_1\otimes\cdots\otimes x_{d-1}\otimes 4tG_d),\\
&\qquad \sum_{I\colon\; I\subset\{1,\dots,d-1\},\, I\neq\emptyset}
\!\!\!   W_I^x(y_d-x_d)|\alpha,4t)
\le \!\! \sum_{I\colon\; I\subset\{1,\dots,d\},\, d\in I,\,
I\cap\{1,\dots,d-1\}\neq\emptyset} \!\!\!\!
W_I^x(\alpha,4t) \biggr\},
\endalign
$$
$$
\align
B&=\biggl\{y=(y_1,\dots,y_d)\colon\; y\in
R^{n_1}\times\cdots\times R^{n_d},\; \\
& \qquad \alpha((y_1\otimes\cdots\otimes y_{d-1}-
x_1\otimes\cdots\otimes x_{d-1})\otimes y_d)
\qquad \!\!\!\! \le \!\!\!\!
\sum_{I\colon\;I\subset\{1,\dots,d-1\},\, I\neq\emptyset} \!\!
W_I^x(y_d|\alpha,4t)\biggr\}
\endalign
$$
and
$$
C=B\cap\{y=(y_1,\dots,y_d)\colon\; y\in
R^{n_1}\times\cdots\times R^{n_d},\; y_d\in A\}.
$$
 
I claim that the inequalities
$$
\mu_{n_d,t}(A)\ge \frac12e^{-1/2t^2} \tag A2
$$
and
$$
\aligned
&\mu_{n_1+\cdots+n_{d-1},t}(B\cap (
(R^{n_1}\times\cdots\times R^{n_{d-1}})\times y_d))
\ge 2^{-(d-1)}e^{-(d-1)/2t^2} \\
&\qquad\qquad \text{for all } y_d=(y_d(1),\dots,y_d(n_d))\in R^{n_d}
\endaligned \tag A3
$$
hold, where $(R^{n_1}\times\cdots\times R^{n_{d-1}})\times y_d=
\{(y_1,\dots,y_{d-1},y_d)\colon\;
(y_1,\dots,y_{d-1})\in R^{n_1}\times\cdots\times R^{n_{d-1}}\}$.

Relation (A2) follows from the identity
$4E\alpha(tG_n)=E\alpha(4tG_n)$ and Proposition~4.1 with the choice
$
\alpha_1(z)=\alpha(x_1\otimes\cdots\otimes x_{d-1}\otimes z)$
and $\alpha_2(z)
=\summ_{I\colon\; I\subset\{1,\cdots,d-1\},\, I\neq\emptyset}
W^x_I(z|\alpha,4t)
$
for $z\in R^{n_d}$. Observe that both $\alpha_1(\cdot)$ and
$\alpha_2(\cdot)$ are pseudonorms in $R^{n_d}$, hence Proposition~(4.1)
is applicable for them.
 
Relation (A3) follows from Proposition~4.2 with parameter $d-1$ if
it is applied for the pseudonorm $\bar\alpha_{y_d}$
on $R^{n_1}\otimes\cdots\otimes R^{n_{d-1}}$ defined by the
formula $\bar\alpha_{y_d}(u)=\alpha(u\otimes y_d)$ for
$u\in R^{n_1}\otimes\cdots\otimes R^{n_{d-1}}$ with a fixed
$y_d\in R^{n_d}$. Here $u\otimes y_d$ is that function in
$R^{n_1}\otimes\cdots\otimes R^{n_d}$ for which
$u\otimes y_d(i_1,\dots,i_d)=u(i_1,\dots,i_{d-1})y_d(i_d)$ for
all $1\le i_j\le n_j$, $1\le j\le d$. Observe that $\bar\alpha_{u_d}$
is really a pseudonorm for all $y_d\in R^{n_d}$, hence we can
apply Proposition~4.2 with parameter $d-1$ for it.
 
Relations (A2), (A3) and the Fubini theorem imply that
$$
\mu_{n_1+\cdots+n_d,t}(C)\ge2^{-d}e^{-d/2t^2}. \tag A4
$$
Indeed, $\mu_{n_1+\cdots+n_{d-1},t}
(B\cap((R^{n_1}\times\cdots\times R^{n_{d-1}})\times y_d))
\ge 2^{-(d-1)}e^{-(d-1)/2t^2}$ for all points $ y_d\in R^{n_d}$
by relation (A3). We get relation (A4) from this inequality, 
relations (A2) and the Fubini theorem by integrating this inequality 
on the set $\{y_d\in A\}$ with respect to the measure~$\mu_{n_d,t}$.
 
Finally, I claim that
$$
C\subset\left\{y=(y_1,\dots,y_d)\colon\;
y\in R^{n_1}\times\cdots\times R^{n_d},\;
\rho_\alpha(x,y)\le\sum_{I\colon\; I\subset\{1,\dots,d\},\, I\neq\emptyset}
W^x_I(\alpha,4t)\right\}. \tag A5
$$
Indeed, if $y=(y_1,\dots,y_d)\in C$, then
$$
\align
\rho_\alpha(x,y)&
\le \alpha(x_1\otimes\cdots\otimes x_{d-1}\otimes(y_d-x_d))+
\alpha((y_1\otimes\cdots\otimes y_{d-1}-x_1\otimes\cdots x_{d-1})\otimes y_d)\\
&\le E\alpha(x_1\otimes\cdots\otimes x_{d-1}\otimes 4tG_d)+
\sum_{I\colon\;I\subset\{1,\dots,d-1\},\, I\neq\emptyset}
W_I^x(y_d|\alpha,4t)\\
&\le E\alpha(x_1\otimes\cdots\otimes x_{d-1}\otimes 4tG_d)+
\sum_{I\colon\;I\subset\{1,\dots,d-1\},\, I\neq\emptyset}
W_I^x(y_d-x_d|\alpha,4t)\\
&\qquad +\sum_{I\colon\;I\subset\{1,\dots,d-1\},\, I\neq\emptyset}
W_I^x(\alpha,4t)\\
&\le E\alpha(x_1\otimes\cdots\otimes x_{d-1}\otimes 4tG_d)+
\sum_{I\colon\; I\subset\{1,\dots,d\},\, d\in I,\,
I\cap\{1,\dots,d-1\}\neq\emptyset}  W_I^x(\alpha,4t) \\
&\qquad +\sum_{I\colon\;I\subset\{1,\dots,d-1\},\, I\neq\emptyset}
W_I^x(\alpha,4t)
=\sum_{I\colon\;I\subset\{1,\dots,d\},\, I\neq\emptyset} W_I^x(\alpha,4t).
\endalign
$$
The first inequality of this series of inequalities holds because of
relation (A1). The second inequality was based on the first relation
in the definition of the set $A$ and on the definition of the
set $B$. The third inequality is valid because of the relation
$\alpha(z\otimes y)\le\alpha(z\otimes x)+\alpha(z\otimes(y-x))$ for
arbitrary pseudonorm $\alpha$ on the tensor product
$R^{n_1}\otimes\cdots\otimes R^{n_d}$ and
$z\in R^{n_1}\otimes\cdots\otimes R^{n_{d-1}}$. 
The last inequality follows from the second relation in the definition
of the set~$A$. In the closing step we have applied the identity
$E\alpha(x_1\otimes\cdots\otimes x_{d-1}\otimes 4tG_d)
=W^x_{\{d\}}(\alpha,4t)$.
 
Proposition~4.2 is a simple consequence of relations~(A4) and~(A5).

\vfill\eject
 
\medskip\noindent
{\bf References:}
 
\item{1.)} Adamczak, R. (2006) Moment inequalities for
$U$-statistics. {\it Annals of Probability} {\bf34}, 2288--2314
\item{2.)} Lata\l{a}, R. (2006) Estimates of moments and tails of
Gaussian chaoses. {\it Annals of Probability} {\bf34} 2315--2331
\item{3.)} Ledoux, M. (2001) The concentration of measure phenomenon.
{\it Mathematical Surveys and Monographs}\/ {\bf 89} American Mathematical
Society, Providence, RI.
\item{4.)} Major, P. (2006) An estimate on the maximum of a nice
class of stochastic integrals. {\it Probability Theory
and Related Fields.} {\bf 134}, 489--537
\item{5.)} Major, P. (2007) On a multivariate version of
Bernstein's inequality {\it Electronic Journal of Probability}\/ 
{\bf 12} 966--988

\bye